\documentstyle[12pt]{article}
\newlength{\abstractwidth}
\renewcommand{\thefootnote}{\fnsymbol{footnote}}
\renewcommand{\thanks}[1]{\footnote{#1}} 
\newcommand{\starttext}{ \setcounter{footnote}{0}
\renewcommand{\thefootnote}{\arabic{footnote}}}

\newcommand{\be}{\begin{equation}}
\newcommand{\bea}{\begin{eqnarray}}
\newcommand{\eea}{\end{eqnarray}} \newcommand{\ee}{\end{equation}}
 \newcommand{\<}{\langle}
\renewcommand{\>}{\rangle} \def\ba{\begin{eqnarray}}
\def\ea{\end{eqnarray}}

\def\o{\omega}
\def\Re{{\rm Re}}
\def\Im{{\rm Im}}

\def\log{\,{\rm log}\,}

\def\o{\omega}

\def\al{\alpha}
\def\b{\beta}
\def\g{\gamma}

\def\o{\omega}

\def\si{\sigma}

\def\na{\nabla}

\def\ge{\geq}
\def\le{\leq}

\def\ov{\overline}

\def\I{\int}

\def\na{{\nabla}}
\def\hna{\hat\nabla}

 \def\v{\vskip .1in}

\def\[{{\bf [}}
\def\]{{\bf ]}}

\def\pl{\partial}

\begin{document}
\starttext \baselineskip=18pt \setcounter{footnote}{0}
\newtheorem{theorem}{Theorem}
\newtheorem{corollary}{Corollary}
\newtheorem{proposition}{Proposition}
\newtheorem{lemma}{Lemma}
\newtheorem{definition}{Definition}
\begin{center}
{\Large \bf THE MODIFIED K\"AHLER-RICCI FLOW AND SOLITONS
\footnote{Research supported in part by National Science Foundation
grants DMS-07-57372, DMS-06-04805, DMS-05-14003, and DMS-08-48193.
 }}
\\
\bigskip
\bigskip

{\large D.H. Phong$^*$, Jian Song$^{**}$, Jacob Sturm$^\dagger$ and
Ben Weinkove$^\ddagger$} \\

\bigskip

\begin{abstract}
{\small We investigate the K\"ahler-Ricci flow modified by a holomorphic vector field.
We find equivalent analytic criteria for the convergence of the flow to a K\"ahler-Ricci soliton.
In addition,   we relate the asymptotic behavior of the scalar curvature along the flow to the
lower boundedness of the modified Mabuchi energy. }

\end{abstract}

\end{center}

\setlength\arraycolsep{2pt}

\baselineskip=15pt
\setcounter{equation}{0}
\setcounter{footnote}{0}

\section{Introduction}
\setcounter{equation}{0}

Let $M$ be a compact K\"ahler manifold of complex dimension $n$ with $c_1(M)>0$.
A K\"ahler-Ricci soliton on $M$ is a K\"ahler metric
$\omega = \frac{i}{2} g_{\bar kj} dz^j \wedge d\ov{z^k}$  in the cohomology class
$\pi\,  c_1(M)$ together with a holomorphic vector field $X$ such that
\be \label{soliton1}
\textrm{Ric}(\omega) - \omega = {\cal L}_X \omega.
\ee
Alternatively, in coordinate notation, writing $X_{\bar k} = g_{ \bar k \ell} X^{\ell}$,
\be
\label{soliton}
R_{\bar kj}-g_{\bar kj}=\na_jX_{\bar k}.
\ee
Let $\Phi_t$ be the 1-parameter group of automorphisms of $M$
generated by $\Re\,X$.
The family of metrics $g_{\bar kj}(t)\equiv\Phi_{-t}^*(g_{\bar kj})$
provides then a solution of the K\"ahler-Ricci flow,
\be
\label{KRflow}
\dot g_{\bar kj}(t)=-R_{\bar kj}+g_{\bar kj}
\ee
where the evolution in time is just by reparametrization.

\smallskip
If $X$ is the zero vector field then (\ref{soliton1}) reduces to the K\"ahler-Einstein equation.
K\"ahler-Ricci solitons are in many ways similar to extremal metrics, which generalize constant scalar
curvature K\"ahler metrics and are  characterized by the condition
 that  the vector field $\nabla^i R$ is holomorphic.
Extremal metrics are  the critical points of the Calabi functional
$C(g_{\bar kj})=\|R-\bar R\|_{L^2}^2$, where $\bar R$ is the average of the scalar curvature.
A classic conjecture of Yau \cite{Y2} asserts that
the existence of constant scalar curvature metrics in a given integral K\"ahler class should
be equivalent to the stability of the polarization in the sense of geometric invariant
theory. Notions of K-stability for constant scalar curvature metrics
have been proposed by Tian \cite{T} and Donaldson
\cite{D1}, and extended to the case of extremal metrics by Szekelyhidi \cite{S1} (see also \cite{M}).
Similarly, the existence of K\"ahler-Ricci solitons is expected to be equivalent to
a suitable notion of stability.

\smallskip
K\"ahler-Ricci solitons can also be viewed as the stationary points of the
{\it modified K\"ahler-Ricci flow}, that is,
\be
\label{modifiedKRflow}
\dot g_{\bar kj}
=
-R_{\bar kj}+g_{\bar kj}+\na_j X_{\bar k}
\ee
which is the flow
(\ref{KRflow}) reparametrized by the automorphisms $\Phi_t$
generated by the real part $\Re \,X$ of the holomorphic vector field $X$.
Similar reparametrizations, in the context of Hamilton's original Ricci flow \cite{H},
had been introduced by DeTurck \cite{DeT} to give a simpler proof of the short-time
existence of the flow.

\smallskip

The modified K\"ahler-Ricci flow appears in the work of Tian-Zhu \cite{TZ2} as part of their
study of the K\"ahler-Ricci flow assuming \emph{a priori} the presence of a K\"ahler-Ricci soliton.
They make use of a Moser-Trudinger type inequality from \cite{CTZ} to deduce convergence of the flow
in the sense of Cheeger-Gromov.  (In the special case where there are no nontrivial holomorphic
vector fields, it is known by the work of \cite{P2},  \cite{TZ2}, \cite{PSSW2}
that the existence of a K\"ahler-Einstein metric implies the exponential convergence
of the K\"ahler-Ricci flow to that metric.)

\medskip

In this paper, we study the long-time behavior of the modified K\"ahler-Ricci flow
without assuming the existence of a K\"ahler-Ricci soliton.   We
give analytic conditions which are both necessary and sufficient for the
convergence of the flow to a K\"ahler-Ricci soliton.   These conditions are analogous to the ones given in
\cite{PSSW1} for the convergence of the K\"ahler-Ricci flow. As explained in
\cite{PS1} and \cite{PSSW1} they can be interpreted as stability conditions in an infinite-dimensional
geometric invariant theory, where the orbits are those of the diffeomorphism group
acting on the space of almost-complex structures.

\smallskip
We provide now a description of our results.  We will assume always that  $M$ is a compact
K\"ahler manifold with $c_1(M)>0$ and $X$ is a holomorphic vector field whose
imaginary part $\Im X$ induces an $S^1$ action on $M$. Note that once a maximal compact
subgroup $G$ of $\textrm{Aut}^0(M)$ is fixed, there is a natural choice of such a vector
field $X$ (for more details see Section \S\ref{mabfut} below).

\smallskip

First, we define the notion of the Hamiltonian $\theta_{X,\o}$ and modified Ricci potential $u_{X,\o}$.
Write ${\cal K}_X$ for the space of K\"ahler metrics in $\pi c_1(M)$ which are invariant under
$\Im X$.
Given $\omega = {i \over 2} g_{\ov{k} j} dz^j \wedge dz^{\ov{k}} \in {\cal K}_X$
we define a real-valued function $\theta_{X, \omega}$ by
\begin
{eqnarray}
\label{theta}
X^j g_{\ov{k} j} =  \partial_{\ov{k}} \theta_{X, \omega}, \ \
\int_M e^{\theta_{X, \omega}} \omega^n   =  \int_M \omega^n =: V.
\end{eqnarray}
The Ricci potential $f=f(\omega)$ is given by
\begin{eqnarray} \label{riccipotential}
g_{\ov{k} j} - R_{\ov{k} j}  =  \partial_{\ov{k}} \partial_j f, \ \
\int_M e^{-f} \omega^n   = V,
\end{eqnarray}
and we define a modified Ricci potential $u_{X,\o}$ by
\be
\label{modifiedRiccipotential}
u_{X,\o}=f+\theta_{X,\o}.
\ee
If $M$ admits a K\"ahler-Ricci soliton $\omega \in \pi \, c_1(M)$ with respect to the vector field $X$,
then $\omega$ is necessarily in ${\cal K}_X$ and $u_{X,\o}=0$.  Let
$g_{\bar kj}(t)$ evolve by the modified K\"ahler-Ricci flow  and set
\be
\label{YX}
Y_X(t)
=
\int_M |\na u_{X, \o}|^2\,e^{\theta_{X,\o}}\,\o^n.
\ee
The modified K\"ahler-Ricci flow starting at $\omega_0 \in {\cal K}_X$ preserves the K\"ahler class, and can be
expressed as a flow of K\"ahler potentials.   Define
\be
{\cal P}_X (M, \omega_0) = \{ \varphi \in C^{\infty}(M) \ | \ \omega =
\omega_0 + \frac{i}{2} \partial \ov{\partial} \varphi >0, \ \Im X (\varphi)=0 \},
\ee
which, modulo constants, can be identified with  ${\cal K}_X$.  Let $\varphi= \varphi(t) \in {\cal P}_X(M, \omega_0)$ be the solution of the equation
\bea
\label{modifiedKRpotential}
\dot\varphi&=&\log{\o^n\over \o_0^n}+\varphi+\theta_{X,\o}+f(\o_0),
\nonumber\\
\varphi(0)&=& c_0.
\eea
Then the K\"ahler metrics $\omega = \omega_0 + \frac{i}{2} \partial \ov{\partial} \varphi$ evolve by the modified K\"ahler-Ricci flow (\ref{modifiedKRflow}).  
The initial constant $c_0$ can affect the growth of $\varphi$ for large time,
and has to be chosen with some care. We choose it to be given by the
value (\ref{c0}) described in Section \S 2 below.

\medskip
Our first main theorem is a general characterization of the
convergence of the modified K\"ahler-Ricci flow, which shows in
particular that if convergence occurs, it always does so
at an exponential rate:

\begin{theorem}
\label{conditions}
Let $\o_0\in {\cal K}_X$,
$\o_0:= {i\over 2}g_{\bar kj}^0dz^j\wedge d\bar z^k$,
and consider the modified K\"ahler-Ricci flow
(\ref{modifiedKRflow}) with initial metric $\omega_0$.
Then the following conditions are equivalent:

\smallskip
{(i)} The modified K\"ahler-Ricci flow $g_{\bar kj}(t)$ converges in $C^\infty$
to a K\"ahler-Ricci soliton $g_{\bar kj}(\infty)$ with respect to $X$.

\smallskip
{(ii)} The function $\|R-n-\na_jX^j\|_{C^0}$ is integrable, i.e.,
\be
\int_0^\infty\|R-n-\na_jX^j\|_{C^0}\,dt<\infty.
\ee

\smallskip
{(iii)}
Let the potential $\varphi(t)$ 
evolve according to (\ref{modifiedKRpotential}),
with initial value $c_0$ as specified in the
equation (\ref{c0}) below.
Then we have
\be
{\rm sup}_{t\geq 0}\|\varphi(t)\|_{C^0}
<\infty.
\ee

\smallskip
{ (iv)} Let the function $Y_X(t)$ be defined by
(\ref{YX}). Then there exist constants $\kappa, C$ with $\kappa$
strictly positive so that
\be
\label{exponentialdecay}
Y_X(t)\,\leq C\,e^{-\kappa \,t}.
\ee

\smallskip
{(v)}
The modified K\"ahler-Ricci flow $g_{\bar kj}(t)$
converges exponentially fast in $C^\infty$ to a K\"ahler-Ricci
soliton $g_{\bar kj}(\infty)$ with respect to $X$.
\end{theorem}

The preceding theorem shows that the growth
of $Y_X(t)$, or alternatively, of the function $\|R-n-\na_jX^j\|_{C^0}(t)$ is key to the
problem of convergence of the modified
K\"ahler-Ricci flow.  Our next result addresses the behavior of these quantities under a stability assumption.
Following \cite{TZ1}, we define the modified Mabuchi K-energy
$\mu_X: {\cal P}_X(M, \omega_0) \rightarrow \mathbf{R}$ by
\be \label{Kenergy0}
\delta \mu_X (\varphi) = - \frac{1}{V} \int_M \delta \varphi \left( R - n - \nabla_j X^j - X u_{X, \omega} \right) e^{\theta_{X, \omega}} \omega^n, \ \ \mu_X(0)=0.
\ee
To clarify this definition, since  $R - n - \nabla_j X^j - Xu =-(\Delta + \Re X) u_{X, \o}$, the integrand is real and thus $\mu_X$ does indeed map into $\mathbf{R}$.  For a proof that $\mu_X(\varphi)$ is independent of choice of  path in ${\cal P}_X(M, \omega_0)$, see \cite{TZ1}. 

\smallskip

We consider the condition:

\bigskip

($A_X$) \ \ $\mu_X$ is bounded from below on ${\cal P}_X(M, \omega_0)$

\bigskip

In \cite{TZ1} it is shown that ($A_X$) is a necessary condition for the existence of a K\"ahler-Ricci soliton $\omega$
with respect to $X$. Here we shall establish the following theorem:

\begin{theorem}
\label{Rconvergence}
Assume that Condition $(A_X)$ holds, and let $\o_0\in {\cal K}_X$.
Then we have, along the modified K\"ahler-Ricci flow (\ref{modifiedKRflow}) starting at $\omega_0$,
\be
Y_X(t) \rightarrow 0 \quad \textrm{and} \quad \| R-n-\na_jX^j \|_{C^0}  \rightarrow 0,
\quad as\ t\to\infty.
\ee
Furthermore, for any $p>2$, we have
\be
\label{pintegrability}
\int_0^\infty\|R-n-\na_jX^j\|_{C^0}^pdt <\infty.
\ee
\end{theorem}

\smallskip

Note that a metric $\omega \in {\mathcal K}_X$ satisfies
$R-n - \nabla_j X^j =0$ if and only if $\omega$ is a K\"ahler-Ricci soliton with respect to $X$.
However, the convergence of $\| R-n-\na_jX^j \|_{C^0}$ to zero is of course weaker than the convergence of
the metrics $g_{\bar kj}(t)$ themselves to a K\"ahler-Ricci soliton, and this
is to be expected, since the condition $(A_X)$ is only a semi-stability type
of condition.

\smallskip

Next, we describe another consequence of the condition $(A_X)$.  Associated to the modified Mabuchi K-energy is the modified Futaki invariant $F_X$
(see \cite{TZ1}), given by
\be
\label{modifiedFutaki}
F_X(Z)
=
-\int_M (Zu_{X,\o}) \,e^{\theta_{X,\o}}\o^n,
\ee
for a holomorphic vector field $Z$.
The modified Futaki invariant $F_X$ is independent of the choice of
$\o\in{\cal K}_X$. It follows immediately that $F_X\equiv 0$ is a necessary
condition for the existence of a K\"ahler-Ricci soliton in ${\cal K}_X$.

\smallskip

In the unmodified case, corresponding to $X=0$,
the condition $(A_X)$ reduces to the condition $(A)$
of lower boundedness of the Mabuchi K-energy. It is then easy
to show that $(A)$ implies that the unmodified Futaki
invariant $F_{ X=0 }(Z)$ vanishes for all
holomorphic vector fields $Z\in H^0(M,T^{1,0})$,
by differentiating the functional along the integral paths
of $Z$.   We show how to rework this argument to  
prove the analogous statement when $X\neq 0$ (to our knowledge, this result is not in the literature).

\begin{proposition}
\label{Proposition}
If $(A_X)$ holds then $F_X (Z)=0$ for all holomorphic vector fields $Z$.
\end{proposition}

Our third theorem shows that $(A_X)$ together with one additional assumption give
  necessary and sufficient conditions for
the convergence of the metrics $g_{\bar kj}(t)$ themselves.
Let $\lambda(t)$ be the first positive eigenvalue of the operator
$-g^{j \ov{k}} \nabla_j \nabla_{\ov{k}}$ acting on smooth $T^{1,0}$ vector fields.  Namely,
\be
\label{lambda}
\lambda(t)={\rm inf}_{V\perp H^0(M,T^{1,0})}
{\|\bar\partial V\|^2\over \|V\|^2},
\ee
where
$H^0(M,T^{1,0})$ is the space of holomorphic vector fields on $M$ and we are using the natural $L^2$ inner product induced by $g_{\ov{k} j}(t)$ on the spaces $T^{1,0}$ and $T^{1,-1}$.
This quantity was first introduced in the context of the K\"ahler-Ricci flow in \cite{PS1}.  Recall the following condition from \cite{PSSW1}:

\bigskip
$(S)$
\ \  ${\rm inf}_{t\geq 0}\,\lambda (t)\,>\,0.$

\bigskip

Then we have:

\begin{theorem}
\label{Fullconvergence}
The modified K\"ahler-Ricci flow (\ref{modifiedKRflow}), starting at an arbitrary metric $\o_0 \in {\cal K}_X$,
converges exponentially fast in $C^{\infty}$ to a K\"ahler-Ricci soliton with respect to the holomorphic vector field $X$
 if and only if
the conditions $(A_X)$ and $(S)$ are satisfied.
\end{theorem}

Since condition (S) is invariant under automorphisms, an  immediate consequence of Theorem
\ref{Fullconvergence} is that convergence modulo automorphisms implies full convergence.  More precisely, suppose that $g_{\bar k j}(t)$ is a solution of the modified  K\"ahler-Ricci flow starting at $\omega_0 \in {\cal K}_X$ and assume there  exists a family of automorphisms $\{ \Psi_t\}_{t \in [0, \infty)} $ such that $\Psi_t^*(g_{\bar k j})$ converges to a K\"ahler-Ricci soliton with respect to $X$.  Then $g_{\bar k j}(t)$ converges exponentially fast to a K\"ahler-Ricci soliton with respect to $X$.

\medskip

Finally we discuss in more detail the behavior of $Y_X(t)$ which, as can be seen from
Theorem \ref{conditions}, is key to the convergence of the K\"ahler-Ricci flow.
The next result provides information on the
growth of $Y_X(t)$ for the completely general modified
K\"ahler-Ricci flow and brings to light the obstructions to the convergence of the flow.

\smallskip
It is convenient to introduce a quantity $\lambda_X$ which is uniformly equivalent to the
eigenvalue $\lambda$ described above (see Lemma \ref{equivalencelemma} below).
 Equip the spaces $T^{1,0}$ and $T^{1,-1}$ with the norms
\bea
\label{thetanorms}
\|V\|_\theta^2&=&\int_M g_{\bar kj}V^j\overline{V^k}\,e^{\theta_{X,\o}}\o^n,
\nonumber\\
\|W\|_\theta^2
&=&
\int_Mg_{\bar kj}W_{\bar p}^j\overline{W_{\bar q}^k}g^{q\bar p}\,e^{\theta_{X,\omega}}\o^n,
\eea
with respect to which they can be completed into Hilbert spaces.
Define the eigenvalue $\lambda_X(t)$ by
\be
\label{lambdaX}
\lambda_X(t)={\rm inf}_{V\perp H^0(M,T^{1,0})}
{\|\bar\partial V\|_\theta^2\over \|V\|_\theta^2},
\ee
where
the notion of perpendicularity is taken with respect to the norm $\|\cdot\|_{\theta}$.
Then we have the following:

\begin{theorem}
\label{differencedifferential}
Consider the modified K\"ahler-Ricci flow
(\ref{modifiedKRflow}) with initial metric $\omega_0 \in {\cal K}_X$.
Then there exist $N$, $\delta_j\geq 0$,
$0\leq j\leq N$, depending only on $n$
and satisfying the condition
\be
{1\over 2}\sum_{j=0}^N\delta_j>1
\ee
and a constant $C>0$,
depending only on the initial metric $\o_0$ so that, for all $t \ge 2N$,
the following difference-differential inequality holds 
\bea
\label{estimate}
\dot Y_X(t)
&\leq& -2\,\lambda_X(t)\,Y_X(t)
- 
2\,\lambda_X(t)\,F_X(\pi (\overline{\na}(u_{X,\o})))
+
C\,\prod_{j=0}^N Y_X(t-2j)^{\delta_j\over 2}.
\eea
Here $\overline{\na} u_{X,\o}$ is the vector field $g^{j\bar k}\pl_{\bar k}u_{X,\o}$,
and $\pi$ is the orthogonal projection, with respect to the norm $\|\cdot\|_\theta$,
of the space $T^{1,0}$ of vector fields
onto its subspace
$H^0(M,T^{1,0})$ of holomorphic vector fields.
\end{theorem}

We stress that the estimate
(\ref{estimate})
holds in full generality, without any stability assumptions.

The main point of the estimate (\ref{estimate}) is that certain
difference-differential inequalities, while weaker than the
standard inequality $\dot Y\leq -\lambda Y$, can still guarantee
the exponential decay of $Y(t)$.
A key example for our purposes is
Lemma \ref{formulationlemma} below. In view
of the estimate (\ref{estimate}) for $\dot Y_X(t)$, we see that the convergence of the
modified K\"ahler-Ricci flow is related to three issues:

\smallskip
(a) The vanishing of the modified Futaki invariant $F_X$;

(b) Whether $Y_X(t)$ tends to $0$ as $t\to\infty$;

(c) The existence of a strictly positive uniform lower bound for $\lambda_X(t)$ (or $\lambda(t)$).

\smallskip

The arguments and viewpoint in this paper run closely in parallel to
those of \cite{PSSW1}, which deals with the case $X=0$.
In what follows, we have emphasized only the
main changes due to the modified K\"ahler-Ricci flow and the holomorphic vector
field $X$. Nevertheless, we have tried to make the discussion reasonably self-contained,
and taken the opportunity to bring out some estimates for both the K\"ahler-Ricci
flow and the modified K\"ahler-Ricci flow which hold in all generality,
independent of any stability assumption. That is particularly the case
for Theorems \ref{conditions} and \ref{differencedifferential} above.

\smallskip

In addition, we mention one result used in the proofs of the main theorems which may be of independent
interest.   Assuming only
that the initial metric $g_{\ov{k} j}$ is invariant under $\textrm{Im} X$, we have (see Proposition 2 below)
\be \label{X2bound}
g_{\ov{k} j} X^j \ov{X^k} = |X|^2 \le C,
\ee
for all time, where $g_{\ov{k} j}(t)$ is a solution of either the K\"ahler-Ricci flow or the modified flow.

\medskip

We give a brief outline of the paper.  In Section \S\ref{mKRF} we give some background on the modified
Mabuchi K-energy and Futaki invariant, prove Proposition \ref{Proposition} and describe the K\"ahler-Ricci and modified K\"ahler-Ricci flows.
In Section \S\ref{EKRF} we prove some estimates for the modified K\"ahler-Ricci flow, which hold in a
general setting.   The proof of the difference-differential inequality,
Theorem \ref{differencedifferential}, is given  in Section \S\ref{PT4}.
Theorems \ref{Rconvergence}, \ref{conditions} and \ref{Fullconvergence}
are then proved in Sections \S\ref{PT2}, \S\ref{PT1} and \S\ref{PT3} respectively.
Finally, in Section \S\ref{remarks} we mention a few further remarks and questions.

\pagebreak[3]
\section{The modified Mabuchi K-energy, Futaki invariant and K\"ahler-Ricci flow}
\label{mKRF}
\setcounter{equation}{0}

\subsection{The modified Mabuchi K-energy and Futaki invariant}
\label{mabfut}

We discuss in this subsection the modified Mabuchi K-energy and
modified Futaki invariant and mention two results of Tian-Zhu.

From the definition of the modified Mabuchi K-energy,
\be
\delta \mu_X (\varphi) = - \frac{1}{V} \int_M \delta \varphi
\left( R - n - \nabla_j X^j - X u_{X, \omega} \right) e^{\theta_{X, \omega}}
\omega^n,  \ \ \mu_X(0)=0,
\ee
we see that its critical points satisfy
\be
(\Delta + X) u_{X, \omega} =0,
\ee
and hence  are K\"ahler-Ricci solitons with respect to $X$.  Note that after integrating by parts
the variation of $\mu_X$ can be rewritten as
\be  \label{Kenergy1}
\delta\mu_X (\varphi) = - \frac{n}{V} \int_M \frac{i}{2}
\partial u_{X,\o} \wedge
\overline{\partial} {\delta\varphi} \wedge e^{\theta_{X, \omega}} \omega^{n-1},
\ee
a formulation that will be useful later.

\medskip

A key result that we will use in this paper is as follows:

\begin{theorem}
(\cite{TZ1})
\label{theoremlowerbound}
If there exists a K\"ahler-Ricci soliton with respect to $X$,
then $\mu_X$ is bounded below on ${\cal P}_X (M, \omega_0)$.
\end{theorem}

Next, the modified Futaki invariant $F_X: H^0(M, T^{1,0}) \rightarrow \mathbf{C}$ is given by
\be
F_X(Z)
=
-\int_M (Zu_{X,\o}) \,e^{\theta_{X,\o}}\o^n.
\ee
An important property of this invariant is that  $F_X$ vanishes on the reductive part of $H^0(M, T^{1,0})$ for a particular choice of $X$.  More precisely, fix a maximal compact subgroup $G$ of the identity component $\textrm{Aut}^0(M)$ of the automorphism group of $M$.  Write $\textrm{Lie}\, (G)^{\mathbf{C}}$ for the complexification of the  Lie algebra of $G$, a reductive Lie subalgebra of the space of holomorphic vector fields.  Then:

\begin{theorem} (\cite{TZ1}) \label{theoremX}
There exists a unique holomorphic vector field $X'$ in  $\emph{Lie}\, (G)^{\mathbf{C}} \subset H^{0}(M, T^{1,0})$ satisfying $\Im X' \in \emph{Lie}\, G$ and
\be
F_{X'} (Z)= 0, \quad \textrm{for all } Z \in \emph{Lie}\, (G)^{\mathbf{C}}.
\ee
Moreover, $X'$ lies in the center of $\emph{Lie}\, (G)^{\mathbf{C}}$.
\end{theorem}

Note that the condition $\Im X' \in \textrm{Lie}\, G$ ensures that the imaginary part of $X'$
induces an $S^1$ action on $M$.   In this paper, we do not explicitly take our holomorphic vector field $X$ to be this choice $X'$.  On the other hand, it
 follows from Proposition \ref{Proposition} and the uniqueness part of
Theorem  \ref{theoremX} that once we fix a maximal compact subgroup $G$,
if condition $(A_X)$ holds and  $\Im X$ is in $\textrm{Lie}\, G$ then $X=X'$.

Theorem \ref{theoremX} is related to the uniqueness result of \cite{TZ1} for K\"ahler-Ricci solitons, which we state now for the reader's convenience.  If $\omega$, $\omega'$ are K\"ahler-Ricci solitons with respect to holomorphic vector fields $X$, $X'$ then there exists $\sigma$ in $\textrm{Aut}^0(M)$ such that
\be
\omega = \sigma^* \omega', \ \ \textrm{and} \ \ X' = \sigma_* (X).
\ee

\subsection{$(A_X)$ and the modified Futaki invariant $F_X$}

In this subsection, we give the proof of Proposition \ref{Proposition}.  

\medskip

We will first show that $F_X(Z)=0$ for all holomorphic vector fields $Z$ satisfying
\be \label{Z}
{\cal L}_{\Im X} Z =0.
\ee
Fix a K\"ahler metric $\omega_0 \in {\cal K}_X$.  Write $\Psi_t$ for the 1-parameter family of automorphisms of $M$ induced by $\Re Z$.
Define $\omega_t$ and $\psi_t$ by
\be
\Psi_t^* \omega_0 =\omega_t = \omega_0 + \frac{i}{2} \partial \bar\partial \psi_t.
\ee
Note that $\psi_t$ is defined only up to the addition of a constant, but this will not affect our calculations.  
Since by assumption $Z$ is invariant under $\Im X$ we see that $\psi_t \in {\cal P}_X(M, \omega_0)$ and
hence $\mu_X(\psi_t)$ is well-defined.  Also,
\be 
\frac{i}{2} \partial \bar\partial \dot{\psi} = {\cal L}_{\Re Z} \omega,
\ee
where here, and henceforth, we are dropping the $t$ subscript.

\smallskip

On the other hand, there exists a complex-valued function
$\theta_{Z, \o}$, invariant under $\Im X$,  such that
\be
\iota_{Z}\omega= \frac{i}{2} \bar\partial \theta_{Z, \o}. \label{dbarZ}
\ee
Indeed, all manifolds with positive first Chern class are simply connected so
the $\bar \partial$-closed $(0,1)$-form $Z^j g_{\ov{k} j} d\ov{z^k}$ must be $\bar \partial$-exact.
Since
\be
\frac{i}{2}\partial\bar\partial \dot{\psi}= d \iota_{\Re Z}
\omega= \Re (d \iota_{Z} \omega)=\frac{i}{2} \partial\bar\partial \Re \, \theta_{Z, \o},
\ee
we can assume without loss of generality that
$\dot{\psi}=\Re\,  \theta_{Z, \o}.$
Calculate
\bea
\nonumber
{\cal L}_{\Re  Z} (e^{\theta_{X, \o}} \omega^n) &=&  d( e^{\theta_{X, \o}} \iota_{\Re Z } \omega^n)\\
\nonumber
&=&e^{\theta_{X, \o}} \left(\Delta \dot{\psi}\, \omega^n + n\,
d\theta_{X, \o}\wedge \iota_{\Re Z}\omega\wedge\omega^{n-1}\right)\\
\label{eqn00}
&=& e^{\theta_{X, \o}} (\Delta \dot{\psi}\,  \omega^n + n\,
\Re (  \frac{i}{2} \, \partial \theta_{X, \o}\wedge \bar\partial\theta_{Z, \o} \wedge\omega^{n-1})).
\eea
Then, by (\ref{Kenergy1}),
\bea \nonumber
\frac{d}{dt} {{\mu}_X} (\psi)
&=& - \Re \left( \frac{n}{V} \int_M  \frac{i}{2}
\partial u_{X, \omega} \wedge\bar\partial \dot{\psi}
\wedge e^{\theta_{X, \o}}\omega^{n-1} \right) \\ \nonumber
&=& \frac{n}{V} \int_M  u_{X, \omega} e^{\theta_{X, \o}} \left( \frac{i}{2}
\partial\bar\partial \dot{\psi} + \Re \left(\frac{i}{2} \partial \theta_{X, \o}
\wedge\bar\partial\dot{\psi} \right) \right)\wedge \omega^{n-1}\\ \nonumber
&=&\frac{1}{V} \int_M u_{X, \o} e^{\theta_{X, \o}} \left( \Delta \dot{\psi} \, \omega^n+ n
\Re \left(\frac{i}{2}\partial \theta_{X, \o}\wedge\bar\partial\theta_{Z, \o}
\wedge \omega^{n-1} \right) \right)\\ \nonumber
&=& \frac{1}{V} \int_M u_{X, \o} {\cal L}_{\Re Z}(e^{\theta_{X, \o}}\omega^n)\\
\nonumber
&=& - \frac{1}{V} \int_M (\Re Z) (u_{X, \o}) e^{\theta_{X, \o}}\omega^n\\
\label{eqn11}
&=&\frac{1}{V} \, \Re \left( F_X(Z) \right).
\eea
To go from the 2nd to 3rd lines, we have used the fact that
$\dot{\psi} = \theta_{Z, \o} - i \Im \, \theta_{Z, \o}$ and
\bea \nonumber
 n \Re ( \frac{1}{2} \partial \theta_{X, \o} \wedge \bar \partial \Im \, \theta_{Z, \o} \wedge\o^{n-1})
& = & - n \Re \left( \frac{1}{2} \partial \Im \, \theta_{Z, \o} \wedge \bar \partial \theta_{X, \o}
\wedge \omega^{n-1} \right) \\ \nonumber
& = &   - (\Im \, X) (\Im \, \theta_{Z, \o}) \, \o^n   =  0,
\eea
since by the assumption (\ref{Z}), $\Im \, \theta_{Z, \o}$ is invariant under $\Im \, X$.

\smallskip

Condition $(A_X)$ implies from (\ref{eqn11}) that
$\Re \left( F_X(Z) \right)=0$. Replacing $Z$ by $i Z$ shows that $F_X(Z)=0$ for all holomorphic vector fields $Z$ invariant under $\Im \, X$.  

\smallskip

Now let $Z$ be an arbitrary holomorphic vector field.
Denote by $\{ \sigma_t \}_{t \in [0, 2\pi]}$ the
1-parameter family of automorphisms
induced by ${\rm Im}\, X$ and define a holomorphic vector field $\hat{Z}$ by
\be \hat Z\ = \ {1\over 2\pi}\I_0^{2\pi}
\sigma^*_t Z\ dt.
\ee
By the argument above, since $\hat{Z}$ is invariant under $\Im \, X$, we have
$F_X(\hat Z)=0$.  Then 
for each  $t \in [0,2\pi]$,
\be F_X(Z) \ = \ -\I_M (Zu_{X,\o_0})\,\o_0^n
\ = \ -\I_M(\si^*_t\,Z)
(\si_t^*u_{X,\o_0}) \,\si_t^*\o_0^n
 \ = \ 
F_X(\si^*_tZ),
\ee
using the fact that $\o_0$ and $u_{X,\o_0}$ are invariant
under ${\rm Im}\,X$. Thus
\be F_X(Z)\ = \ {1\over 2\pi}\I_0^{2\pi}
F_X(\si^*_tZ)\ dt\ = \ F_X(\hat Z)\ = \ 0,
\ee
completing the proof of Proposition \ref{Proposition}.

\subsection{The modified K\"ahler-Ricci flow}

We list now some basic properties of the modified K\"ahler-Ricci flow.
Some result from
the fact that the modified K\"ahler-Ricci flow
is a reparametrization of the K\"ahler-Ricci flow,
and the most important is a normalization for the modified K\"ahler
potential $\varphi$ which ensures that ${\rm sup}_{t\geq 0}\|\dot\varphi\|_{C^0}<\infty$.  We remark that this bound is proved in \cite{TZ2} assuming the existence of a K\"ahler-Ricci soliton.  In our case, we only require the invariance of the initial metric $\omega_0$ under $\Im \, X$.

\medskip

Fix $\omega_0 \in {\cal K}_X$, $\omega_0 = \frac{i}{2} g^0_{\bar k j} dz^j \wedge d\ov{z^k}$.  Then the K\"ahler-Ricci flow and modified K\"ahler-Ricci flow starting at $\omega_0$ are given by
\bea
\label{KRF1}
\frac{\partial}{\partial t} \tilde g_{\bar kj} (t) & = & - \tilde R_{\bar kj}
+
\tilde g_{\bar kj}, \quad \tilde{g}_{\bar k j}(0) = g^0_{\bar k j} \\
\frac{\partial}{\partial t} g_{\bar kj}(t) & = & -R_{\bar kj}+
g_{\bar kj}+\na_jX_{\bar k}, \quad g_{\bar k j}(0) = g^0_{\bar k j},
\label{KRF2}
\eea
respectively.   
Note that if $\{ \Phi_t\}_{t \in [0, \infty)}$, $\Phi_0 =\textrm{id}$, is the 1-parameter family of automorphisms of $M$ generated by $\Re X$,  
then the solutions to (\ref{KRF1}) and (\ref{KRF2}) are related by $g_{\bar kj}(t)= \Phi_t^* ( \tilde{g}_{\bar kj})$.  The K\"ahler-Ricci flow preserves the $S^1$ action induced by $\Im \, X$ and so the K\"ahler
forms $\tilde{\o}(t)$ and $\o(t)$ lie in ${\cal K}_X$.
In the sequel, we will often drop the $t$.   Also, we will denote by $\tilde{f}$, $\tilde{\nabla}$ and $\tilde{\Delta}$ the Ricci potential (see (\ref{riccipotential})), covariant derivative and Laplacian with respect to $\tilde{g}_{\bar k j}$.
Since we are using different notation for solutions of (\ref{KRF1}) and (\ref{KRF2}),
and each of these flows can easily be obtained from the other,  we will sometimes refer simply to
`the flow' rather than the specific equation (\ref{KRF1}) or (\ref{KRF2}).

\smallskip

Before continuing our discussion of these two flows, we list without proof some
known estimates for the K\"ahler-Ricci flow.  The first three statements are due to
Perelman \cite{P2} (see \cite{ST}) and the fourth is due to Zhang \cite{Zha} and Ye \cite{Ye}.  

\begin{theorem}  \label{perelman} In the following, all norms are taken with respect to the metric $\tilde{g}_{\bar k j}(t)$.
\begin{enumerate}
\item[(i)] There exists a constant $C$ depending only on $\omega_0$ such that the Ricci
potential $\tilde{f}$ along the flow satisfies
\be \label{p1}
\| \tilde{f} \|_{C^0} + \| \tilde{\nabla} \tilde{f} \|_{C^0} + \| \tilde{\Delta} \tilde{f} \|_{C^0} \le C.
\ee
\item[(ii)]  The diameters $\emph{diam}_{\tilde{g}(t)} M$ are uniformly bounded  along
the flow by a constant depending only on $\omega_0$.
\item[(iii)]  Let $\rho>0$ be fixed.
Then there exists a constant $c>0$ depending only on $\omega_0$ and $\rho$
such that for all $x \in M$ and all $r$ with $0 < r \le \rho$ we have
\be
\int_{B_r(x)} \tilde{\omega}^n(t) \ge c \, r^{2n},
\ee
where $B_r(x)$ is the geodesic ball centered at $x$ of radius $r$
with respect to $\tilde{g}_{\ov{k} j}(t)$.
\item[(iv)]  There exists a constant  $C_S$ depending only on $\omega_0$ such that the Sobolev inequality
\be
\| \eta \|_{L^{2n/(n-1)}} \le C_S ( \| \tilde{\nabla} \eta \|_{L^2} + \| \eta \|_{L^2}), \quad \textrm{for all } \eta \in C^{\infty}(M)
\ee
holds.
 \end{enumerate}
\end{theorem}

We remark that this theorem makes no reference to the vector field $X$ and indeed
does not require the initial metric $\omega_0$ to be invariant under $\Im X$.
Moreover, all of the above statements  are invariant under automorphisms
and hence the analogous statements hold also for the metrics $g_{\ov{k} j}$.

\smallskip

We now describe (\ref{KRF1}) and (\ref{KRF2}) in terms of potentials.
Define $\tilde{\varphi} = \tilde{\varphi}(t)$ and $\varphi=\varphi(t)$ in ${\cal P}_X (M, \omega_0)$ by
\bea \label{KRFP1}
\frac{\partial \tilde{\varphi}}{\partial t}   & = &  \log \frac{\tilde{\omega}^n}{\omega_0^n} + \tilde{\varphi} + f(\omega_0),  \quad  \tilde{\varphi}(0)= \tilde{c}_0, \\
\frac{\partial \varphi}{\partial t}    & = &  \log \frac{\omega^n}{\omega_0^n} + \varphi + \theta_{X, \omega} + f(\omega_0),  \quad \varphi(0) = c_0, \label{KRFP2}
\eea
where the constants $\tilde{c}_0$ and $c_0$ will be defined shortly.
One can check that
$\tilde{\omega} = \omega_0 +  \frac{i}{2} \partial \ov{\partial} \tilde{\varphi}$ and
$\omega = \omega_0 + \frac{i}{2} \partial \ov{\partial} \varphi$
satisfy (\ref{KRF1}) and (\ref{KRF2}) respectively.

\smallskip
Recall that  $\theta_{X, \omega}$ is the Hamiltonian function defined by (\ref{theta}).
It is well-known that  $\theta_{X,\omega}$ is well-defined.   Indeed, by the same argument as for (\ref{dbarZ}), there is a complex-valued function $\theta_{X, \omega}$ solving
$X^j g_{\bar k j} = \partial_{\bar k} \theta_{X, \omega}$.
The equalities
 \be
 0 = {\cal L}_{\Im X} \omega = d \, \iota_{ \Im  X} \omega =
\frac{i}{2} \partial \bar \partial \,  \Im \, \theta_{X, \omega}
\ee
ensure the existence of a real-valued function  $\theta_{X, \omega}$ satisfying (\ref{theta}).

\smallskip
 Following the conventions of  \cite{PSS} (cf. \cite{CT})
we define
\be
\tilde{c}_0: =  \frac{1}{V} \int_0^{\infty} e^{-t} \int_M | \tilde{\nabla}
\tilde{f} |^2  \tilde{\omega}^n dt - \frac{1}{V} \int_M f(\omega_0) \omega_0^n.
\ee
Note that by Theorem \ref{perelman} we have
$| \tilde{\nabla}  \tilde{f}| \le C$ and so $\tilde{c}_0$ is finite.
With this choice of $\tilde{c}_0$, there exists a uniform $C$ such that
\be \label{ddtphi}
\| \partial \tilde{\varphi}/ \partial t \|_{C^0} \le C.
\ee

Before we define $c_0$, we will need two lemmas.  
\smallskip
The first lemma implies that $\theta_{X, \omega}$ is uniformly bounded along the flow.

\begin{lemma} \label{lemmathetaX}
For all $\omega' \in {\cal K}_X$, we have
\begin{equation}
\left\| \theta_{X, \omega'}    \right\|_{C^0} =  \left\| \theta_{X, \omega_0}    \right\|_{C^0}.
\end{equation}
\end{lemma}
{\it Proof of Lemma \ref{lemmathetaX}:} \ This result is well-known (see \cite{FM} or \cite{Zhu1}, for example), but for the reader's convenience we outline here a proof.  We use the following version of  Moser's theorem  (which can be easily derived from \cite{CdS}, p.43-44, for example).

\begin{theorem} \label{moser} Let $M$ be a compact complex manifold with K\"ahler forms $\omega_0$ and $\omega_1$ which are invariant under an $S^1$ action induced by a real vector field $W$.  Assume in addition that $[\omega_0]=[\omega_1]$.  Then there exists a diffeomorphism $\Psi$ of $M$ satisfying
\be
\Psi^* \omega_1= \omega_0 \quad \textrm{and} \quad \Psi_* W=W.
\ee
\end{theorem}

Write $\omega_1 := \omega' \in {\cal K}_X$.  By definition, $\theta_0 := \theta_{X, \omega_0}$ and $\theta_1 := \theta_{X, \omega_1}$ satisfy
\be
\iota_{\Im X} \omega_0  = \frac{1}{4} d \theta_0, \quad \iota_{\Im X} \omega_1 = \frac{1}{4} d\theta_1,
\ee
with
\be \label{normal}
\int_M e^{\theta_0} \omega_0^n = V = \int_M e^{\theta_1} \omega_1^n.
\ee
From Theorem \ref{moser} we obtain a diffeomorphism $\Psi$ of $M$ with
$\Psi^* \omega_1= \omega_0$ and $\Psi_*{\Im X} = \Im X$.
Applying $\Psi^*$  to the equation $\iota_{\Im X} \omega_1 = \frac{1}{4} d\theta_1$ we obtain
$d \Psi^* \theta_1 = d\theta_0,$
and hence
$\Psi^*\theta_1 = \theta_0 + c$
for some constant $c$.  But from (\ref{normal}) we see that $c=0$.
This implies that $\theta_0, \theta_1: M \rightarrow \mathbf{R}$ have the same image in $\mathbf{R}$.  Q.E.D.

\bigskip

Given this lemma we can now prove the following.

\begin{lemma}
\label{lemmaX}
Along the modified K\"ahler-Ricci flow we have
\be
\int_M |X|^2 e^{\theta_{X, \omega}} \omega^n \le C.
\ee
\end{lemma}
{\it Proof of Lemma \ref{lemmaX}:} \ From the definition of the modified Futaki invariant  and the definition of $\theta_{X, \omega}$ we see that
\be
\int_M |X|^2 e^{\theta_{X, \omega}} \omega^n = - \int_M (Xf) e^{\theta_{X, \omega}} \omega^n - F_X(X).
\ee
Hence, since $F_X(X)$ is independent of choice of metric, we have
\begin{eqnarray} \nonumber
\int_M |X|^2 e^{\theta_{X, \omega}} \omega^n & \le & \int_M |X| | \nabla f| e^{\theta_{X, \omega}} \omega^n + C \\
& \le & \frac{1}{2} \int_M |X|^2 e^{\theta_{X, \omega}} \omega^n + \frac{1}{2} \int_M | \nabla f|^2 e^{\theta_{X, \omega}} \omega^n + C,
\end{eqnarray}
and the lemma follows from Theorem \ref{perelman}, part (i) 
and Lemma  \ref{lemmathetaX}. Q.E.D.

\bigskip

We now define the constant $c_0$ as follows:
\be
\label{c0}
c_0: = \frac{1}{V} \int_0^{\infty} e^{-t} \int_M | \nabla u_{X, \omega}|^2
e^{\theta_{X, \omega}} \omega^n \, dt - \frac{1}{V} \int_M u_{X, \omega_0}
e^{\theta_{X, \omega_0}} \omega_0^n,
\ee
where we recall that $u_{X, \omega}  = f + \theta_{X, \omega}$.
To see that $c_0$ is finite, observe that
\be
| \nabla u_{X, \omega} |^2 \le 2 ( | \nabla f |^2 + |X|^2 ) \le C + 2 |X|^2,
\ee
and hence by Lemma \ref{lemmathetaX} and Lemma \ref{lemmaX},
\be \label{L2u}
\int_M | \nabla u_{X, \omega} |^2 e^{\theta_{X, \omega}} \omega^n \le C.
\ee

We will end this section by proving a uniform bound on  $\dot{\varphi}$.  First, we have another general lemma on Hamiltonian functions:

\begin{lemma}  \label{lemmaTZ}
For any $\omega' \in {\cal K}_X$ with $\omega' = \omega_0 + \frac{i}{2} \partial \bar \partial \varphi'$, the Hamiltonian functions $\theta_{X, \omega'}$ and $\theta_{X, \omega_0}$ are related by:
\be \label{ddttheta}
\theta_{X, \omega'} = \theta_{X, \omega_0}  + X (\varphi'),
\ee
\end{lemma}
{\it Proof of Lemma \ref{lemmaTZ}:} \ This is proved in \cite{TZ1}, p.301. Q.E.D.

\bigskip

We can now prove:

\begin{lemma} \label{lemmaddtphi}
There exists a uniform constant $C$ such that along the flow,
\be \label{ddtphibound}
\| \dot{\varphi}  \|_{C^0} \le C.
\ee
\end{lemma}
{\it Proof of Lemma \ref{lemmaddtphi}:} \ 
From (\ref{KRF1}) and (\ref{KRF2}) we obtain
\be
\label{phis}
\frac{\partial \varphi}{\partial t} = \Phi_t^* \frac{\partial \tilde{\varphi}}{\partial t} +
\theta_{X, \omega} + m(t),
\ee
for some constant $m(t)$.
Define 
\be 
\alpha(t) = \frac{1}{V} \int_M\dot{\varphi}
e^{\theta_{X, \omega}} \omega^n.
\ee
Using Lemma \ref{lemmaTZ}, we have
\be \label{ev1}
\frac{\partial \dot{\varphi}}{\partial t} = (\Delta + X)\dot{\varphi} + \dot{\varphi}.
\ee
Since
\be \label{phiu}
\dot{\varphi} = u_{X, \omega} + c,
\ee
for a constant $c$ depending only on time, we have
\be
\frac{d}{dt}\alpha = \alpha -  \frac{1}{V} \int_M | \nabla u_{X, \o} |^2 e^{\theta_{X, \omega}} \omega^n.
\ee
Integrating this ODE (cf. the argument in \cite{PSS}) and applying (\ref{L2u})  shows that
\be \label{alpha}
0 \le \alpha(t) = \frac{1}{V} \int_t^{\infty} e^{-(s-t)} \int_M | \nabla u_{X, \o} |^2(s)
e^{\theta_{X, \omega(s)}} \omega^n(s) \, ds \le C,
\ee
 or in other words, the average of  $\dot{\varphi}$ with respect to the
measure $e^{\theta_{X, \omega}} \omega^n$ is bounded along the flow.
The lemma now follows immediately from the
formula (\ref{phis}) together with (\ref{ddtphi}), (\ref{alpha}) and Lemma \ref{lemmathetaX}.  Q.E.D.
\bigskip

\pagebreak[3]
\section{Estimates for the modified K\"ahler-Ricci flow} \label{EKRF}
\setcounter{equation}{0}

In this section, we establish the key estimates for the modified K\"ahler-Ricci
flow needed in the sequel. They include the analogues of Perelman's estimates
for the Ricci potential and the scalar curvature for the K\"ahler-Ricci
flow, the estimates for the Laplacian of the Hamiltonian function $\theta_{X, \omega}$, and
a smoothing lemma.

\subsection{Estimates for the modified Ricci potential}

Define
\be
v:= -\dot\varphi.
\ee
Recall from (\ref{phiu}) that $\dot{\varphi}$ differs from $f+\theta_{X,\o}=u_{X,\o}$ by a constant
depending on time. Thus, when computing time evolutions, we have
to distinguish between $-v$ and $u_{X,\o}$, but the difference disappears
whenever $\na$ is applied.  In the last section we established the bound
\be
 \| v \|_{C^0} \le C,
\ee
and in this section we will prove the following further estimates for the modified K\"ahler-Ricci flow.

\begin{proposition} \label{prop2}
Along the flow, the quantities
\be
\| \nabla u_{X,\o} \|_{C^0},  \ \| \Delta u_{X, \o} \|_{C^0}, \ \| X\|_{C^0}, \ \emph{and} \, \  \| \Delta \theta_{X, \omega} \|_{C^0}
\ee
are uniformly bounded by a constant depending only on the initial data.  Here, all norms, covariant derivatives and Laplacians are taken with respect to the evolving metric $g_{\bar k j}(t)$.
\end{proposition}

These bounds will be obtained using the maximum principle.  The proof of this proposition is contained in the following three lemmas.

\begin{lemma} \label{evolution}
We have the following identities along the modified K\"ahler-Ricci flow:
\begin{enumerate}
\item[(i)] $\displaystyle{\frac{\partial v}{\partial t} = ( \Delta +X) v  + v}$
\item[(ii)] $\displaystyle{\frac{\partial}{\partial t} | \nabla v |^2 = (\Delta + X) | \nabla v |^2 -
| \nabla \nabla v |^2 - | \nabla \ov{\nabla} v |^2 + | \nabla v|^2.}$
\item[(iii)] $\displaystyle{\frac{\partial}{\partial t} ( \Delta + X) v = (\Delta + X) (\Delta + X) v +
(\Delta + X) v + | \nabla \ov{\nabla} v |^2.}$
\end{enumerate}
\end{lemma}

\noindent
{\it Proof of Lemma \ref{evolution}:}  This is a straightforward calculation (cf. \cite{CTZ}). Q.E.D.

\medskip

Note that although $X$ is not a real operator on functions on $M$,
the quantities $v$, $| \nabla v |^2$ and $\Delta v$ above are all
invariant under $\Im X$, and so $\Delta +X$ can be replaced by the real operator $\Delta + \Re X$.

\medskip

Using these evolution equations we can give a proof of the following lemma.

\begin{lemma} \label{nablav}
There exists a uniform constant $C$ depending
only on the initial data such that along the flow,
\be
  \| \nabla v \|_{C^0} \le C.
\ee
\end{lemma}
{\it Proof of Lemma \ref{nablav}: } \  This is a straightforward modification of Perelman's maximum principle argument for the bound of the gradient of the Ricci potential (see \cite{ST}, Proposition 6).  Since $v$ is uniformly bounded along the flow by Lemma \ref{lemmaddtphi}, we may choose a constant $B$ such that $v+ B \ge 0$.  Define
$$H=\frac{ |\nabla v|^2}{v+2B}.$$
Compute, using Lemma \ref{evolution},
\bea \label{evolveH}
(\Delta + X - \partial_t) H & = &   \frac{H(H-2B)}{v+2B} -  \frac{2 \Re \left( g^{j \bar k} \partial_j H \partial_{\bar k} v \right)}{v+2B} + \frac{ | \nabla \nabla v|^2 + | \nabla \ov{\nabla} v|^2}{v+2B}.
\eea
Fix $T>0$.  At a maximum point of $H(x,t)$ for $(x,t) \in M \times (0,T]$, the middle term on the right side of (\ref{evolveH}) vanishes and the left hand side of (\ref{evolveH}) is nonpositive.  Then
\be 
\sup_{(x,t) \in M \times [0,T]} H(x,t) \le \max(2B, \sup_{x\in M} H(x,0)),
\ee 
and since $v$ is uniformly bounded, the result follows. Q.E.D.

\bigskip

We can now prove:

\begin{lemma} \label{Xbound}
There exists a uniform constant $C$ depending
only on the initial data such that along the flow,
\be
 \| \nabla \theta_{X, \omega} \|_{C^0} = \| X \|_{C^0} \le C.
 \ee

\end{lemma}
{\it Proof of Lemma \ref{Xbound}:} \ Since
$$ \nabla \Phi_t^* \frac{\partial \tilde{\varphi}}{\partial t}  =   \nabla f$$
and $| \nabla f|$ is bounded by Theorem \ref{perelman}, part (i), the result follows  from (\ref{phis})   and Lemma \ref{nablav}.  Q.E.D.

\bigskip

\begin{lemma} \label{Delta}
There exists a uniform constant $C$ such that along the flow we have
\be
  \| \Delta v \|_{C^0} \le C, \qquad \| \Delta \theta_{X, \omega} \|_{C^0} \le C.
 \ee
\end{lemma}
{\it Proof of Lemma \ref{Delta}:} First note that $|Xv| \leq |X| |\nabla v| \leq C$ by Lemma \ref{nablav} and
Lemma \ref{Xbound}. From Lemma \ref{evolution}, part (iii) we have
\bea \label{minv}
(\Delta + X - \partial_t)((\Delta +X)v) & = & - \Delta v - Xv - | \nabla \bar \nabla v|^2
 \le  - (\Delta v) \left(1 + \frac{\Delta v}{n}\right) + C,
\eea
where we have used the elementary inequality $| \Delta v |^2 \le n| \nabla \overline{\nabla} v|^2$.  Fix an arbitrary $T>0$. At a minimum point of $(\Delta +X)v$ on $M \times (0,T]$ the left hand side of (\ref{minv}) is nonnegative and hence $\Delta v$ is bounded uniformly from below at this point.  This gives the lower bound of $(\Delta +X)v$ along the flow, depending only on the initial data.

To estimate $\| \Delta v\|_{C^0}$, it suffices to prove a uniform upper bound for $(\Delta +X) v$. This argument is similar
to Perelman's estimate of the scalar curvature (see \cite{ST}).  Define
\be
G= \frac{ (\Delta +X) v + 2 |\nabla v|^2 }{ v+2B}
\ee
where $B$ is chosen as in the proof of Lemma \ref{nablav}. Compute
\begin{eqnarray}
 \left( \Delta + X - \partial_t \right) G
&=& - 2 \Re \left( \frac{ \nabla G \cdot \overline{\nabla} v }{ v+2B} \right) +
\frac{  |\nabla\overline{\nabla}v|^2 + 2|\nabla\nabla v|^2}{v+2B}  - \frac{2BG}{(v+2B)}.
\end{eqnarray}
Since $1/(v+2B)$, $ |X v|$ and $ |\nabla v|$ are uniformly bounded, we have
\be
\left( \Delta + X - \partial_t \right) G \geq  -2 \Re \left( \frac{ \nabla G \cdot \overline{\nabla} v }
{ v+2B} \right) +  C_1 |\nabla\overline{\nabla}v|^2 - C_2 | \Delta v| + C_3,
\ee
 for uniform constants $C_1$, $C_2$, $C_3>0$ with $C_1$ uniformly bounded from below away from $0$.
By the maximum principle and a similar argument to the one above, we have
 $(\Delta +X) v \leq C$ for some uniform constant $C$.
This gives the estimate for $\Delta v$. Notice that
\be
\Delta (v+ \theta_{X, \omega})= - \Delta f,
\ee
which is uniformly bounded by Theorem \ref{perelman}, part (i).  It follows that  $\Delta \theta_{X, \omega}$ is uniformly bounded.  Q.E.D.

\bigskip

\subsection{An $L^2/C^0$ Poincar\'e inequality}

Recall that we have
the following Poincar\'e-type inequality on K\"ahler manifolds $(M, \omega)$ with $\omega$ in  $\pi c_1(M)$ (see \cite{F}, or Lemma 2 of \cite{PSSW1})
\be \label{poincare}
\frac{1}{V} \int_M \eta^2 e^{-f} \omega^n \le \frac{1}{V} \int | \nabla \eta |^2 e^{-f} \omega^n,
\ee
for all $\eta \in C^{\infty}(M)$  with $\int_M \eta e^{-f} \omega^n=0$. Define $b=b(t)$ by
\be \label{defnb}
b = \frac{1}{V} \int_M u_{X, \omega} e^{-f} \omega^n.
\ee

Making use of (\ref{poincare}) together with Theorem \ref{perelman}, parts (i) and (iii) we can
prove the following:

\begin{lemma}
\label{l3}
There exists a uniform constant $C$ such that
\be
\|u_{X,\o}-b\|_{C^0}^{n+1}
\,\leq C\,\|\na u_{X,\o} \|_{L^2}\,\|\na u_{X,\o}\|_{C^0}^n.
\ee
\end{lemma}
{\it Proof of Lemma \ref{l3}: } \ See Lemma 3 in \cite{PSSW1}.  Q.E.D.

\subsection{A smoothing lemma}

The following is an analogue of the smoothing lemma from \cite{PSSW1} (see \cite{B}, \cite{CTZ} for  related results.)

\begin{lemma}
\label{smoothinglemma}
There exist positive constants $\delta$ and $K$ depending only on $n$ and the constant $C_X = \sup_{t \in [0, \infty)} \|X\|_{C^0}(t)$ with the following property.
For any $\varepsilon$ with $0< \varepsilon \le \delta$ and any $t_0 \ge 0$, if
\be \label{sl1}
 \| (u_{X,\o}-b)(t_0) \|_{C^0} \le \varepsilon,
\ee
then
\be \label{sl2}
\| \nabla u_{X,\o} (t_0 +2) \|_{C^0} + \| (\Delta + X) u_{X,\o} (t_0+2) \|_{C^0} \le K \varepsilon.
\ee
\end{lemma}
{\it Proof of Lemma \ref{smoothinglemma}: } \ The modified Ricci potential $u_{X, \o}$ evolves by
\be \label{evolveu}
\dot{u}_{X, \o} = (\Delta +X) u_{X,\o} + u_{X, \o} - b(t).
\ee
Indeed,  observe that  $u_{X, \o}$ differs from $-v$ by a constant depending only on time, so applying part (i) of Lemma \ref{evolution} gives (\ref{evolveu}) modulo a constant.  To determine the constant, differentiate the equation $\int_M e^{-f} \o^n = V$ in time and use the relation $\dot{f} = \dot{u}_{X, \o} - Xu_{X, \o}$.

Define constants $c=c(t)$ by
\be
\dot{c} = \dot{b} + c , \ \ c(t_0)=0,
\ee
so that $w= - u_{X,\o} +b - c$ satisfies
\be \label{evolvew}
\frac{\partial w}{\partial t} = (\Delta + X )w + w.
\ee
Moreover, $\| w(t_0) \|_{C^0} \le \varepsilon$.

Now since $v$ and $w$ differ only by a constant, we have from
 Lemma \ref{evolution} that
\be \displaystyle{\frac{\partial}{\partial t} | \nabla w |^2 =
(\Delta + X) | \nabla w |^2 - | \nabla \nabla w |^2 - | \nabla \ov{\nabla} w |^2 + | \nabla w|^2,}
\ee
and
\be
\displaystyle{\frac{\partial}{\partial t} ( \Delta + X) w =
(\Delta + X) (\Delta + X) w + (\Delta + X) w + | \nabla \ov{\nabla} w |^2.}
\ee
We can then modify the argument of Lemma 1 in \cite{PSSW1} to  obtain the result.  Assume without loss of generality that $t_0=0$.  By (\ref{evolvew}), the maximum principle gives $\| w (t) \|_{C^0} \le e^2 \varepsilon$ for $t \in [0,2]$.  From
\be
\frac{\partial}{\partial t} ( e^{-2t} (w^2 + t | \nabla w|^2) \le (\Delta +X) (e^{-2t}(w^2 + t | \nabla w|^2),
\ee
we have $\| \nabla w\|^2_{C^0}(t) \le e^4 \varepsilon^2$ and $\|Xw\|_{C^0}(t) \le C_X e^2 \varepsilon$ for $t \in [1,2]$.  Define
\be
L= e^{-(t-1)} ( | \nabla w|^2 - \varepsilon n^{-1} (t-1) (\Delta +X) w),
\ee
and estimate, using the inequality $(\Delta w)^2 \le n | \nabla \bar \nabla w|^2$,
\be  \label{evolveL}
\frac{\partial}{\partial t} L \le (\Delta +X)L + e^{-(t-1)} n^{-1} ( - \Delta w)( \varepsilon + \Delta w) - e^{-(t-1)}\varepsilon n^{-1} Xw.
\ee
We claim that $L < 2(1+C_X) e^4 \varepsilon^2$ for $t \in [1,2]$.  If the claim is false then, at some point $(x', t') \in M \times (1,2]$, when the inequality first fails, we have $-(\Delta +X)w \ge (1+2C_X) e^4 \varepsilon$.  In particular, at this point, $\varepsilon +  \Delta w \le - \varepsilon$.  It then follows from  (\ref{evolveL}) that $-\Delta w \le | X w| \le C_X e^2 \varepsilon$, a contradiction, thus proving the claim.  Hence, at $t=2$, the inequality
\be
(\Delta +X) w > - 2n (1+C_X) e^5 \varepsilon,
\ee
holds on all of $M$.  To obtain the upper bound of $(\Delta +X)w$ we evolve the quantity
\be
U = e^{-(t-1)} ( | \nabla w|^2 + \varepsilon n^{-1} (t-1) (\Delta +X) w),
\ee
and conclude by a similar argument, after choosing $\varepsilon$ sufficiently small, that $(\Delta +X)w < 2n(1+C_X)  e^5 \varepsilon$ at $t=2$.  Q.E.D.

\bigskip

\noindent
{\bf Remark} \,   Note that by the uniform bound of $|X|$ along the flow,
we can replace (\ref{sl2}) in the conclusion of Lemma \ref{smoothinglemma} with
\be \label{sl3}
\| \nabla u_{X,\o} (t_0 +2) \|_{C^0} + \| \Delta u_{X,\o} (t_0+2) \|_{C^0} \le K \varepsilon.
\ee

\pagebreak[3]
\section{Proof of Theorem \ref{differencedifferential}} \label{PT4}
\setcounter{equation}{0}

We begin by deriving the following analogue for the modified
K\"ahler-Ricci flow of an identity in \cite{PS1},
\bea
\label{analogue1}
\dot Y_X
&=&
- 2\|\bar\na\bar\na u\|_\theta^2
+
\int_M (Xu)|\na u|^2 e^\theta\o^n
\nonumber\\
&&
-
\int_M (R_{\bar kj}-g_{\bar kj}-\na_jX_{\bar k})\na^ju \na^{\ov{k}} ue^\theta
\o^n
-
\int_M(R-n-\na_jX^j)|\na u|^2e^\theta\o^n.
\eea
Here and henceforth, we have denoted $u_{X,\o}$ and $\theta_{X,\o}$ just by $u$ and $\theta$ to
simplify the notation. Note that, as in (\ref{thetanorms}),
the notation $\|\cdot\|_\theta$ refers to $L^2$ norms
with the volume form $\o^n$ replaced by $e^{\theta}\o^n$.
To establish the above identity,
we use the flow for $|\na u|^2= | \nabla v|^2$ already derived in
Lemma \ref{evolution}, from which it follows that
\bea
\label{identity}
\pl_t \left( \int_M|\na u|^2 e^\theta\o^n \right)
&=&
\int_M(\Delta+X)|\na u|^2e^\theta\o^n
-
\int_M|\na\na u|^2e^\theta\o^n-\int_M| \na \bar\na u|^2e^\theta\o^n
\nonumber\\
&&
+\int_M|\na u|^2 e^\theta\o^n
+
\int_M|\na u|^2 (Xu) e^\theta\o^n
\nonumber\\
&&
-
\int_M|\na u|^2(R-n-\na_jX^j)e^\theta\o^n.
\eea
This formula can be simplified as follows. The first term actually vanishes since
\be
\label{zero}
\int_M  e^{\theta}\o^n \,(\Delta+X)\eta =0
\ee
for any smooth function $\eta$,
as can be easily checked by integration by parts, using the relation 
$X_{\bar k}=\pl_{\bar k}\theta$.
Next, we have the following formula of Bochner-Kodaira type,
if $X^j$ is a holomorphic vector field and $u$ is a function invariant
under $\Im X$,
\bea
\label{bk}
\|\nabla\bar\nabla u\|_\theta^2
&=&
\|\bar\nabla\bar\nabla u\|_\theta^2
+
\int_M R_{\bar kj}\nabla^ju\nabla^{\bar k}u\, e^\theta \o^n
-
\int_M \nabla_j X_{\bar k}   \nabla^j u\nabla^{\bar k}u\, e^\theta\o^n.
\eea
To establish this, we integrate by parts,
\bea
\|\nabla\bar\nabla u\|_\theta^2
&=&
\int_M \nabla_j\nabla_{\bar k}u\,\nabla_{\bar p}\nabla_q u\,
g^{j\bar p}g^{q\bar k}\,e^{\theta}\o^n
\nonumber\\
&=&
-
\int_M(\nabla_{\bar p}\nabla_j\nabla_{\bar k}u+X_{\bar p}\nabla_j\nabla_{\bar k}u)\nabla_qu\,
g^{j\bar p}g^{q\bar k}\,e^{\theta}\o^n
\nonumber\\
&=&
-
\int_M\bigg(-R_{\bar pj\bar k}{}^{\bar \ell}\nabla_{\bar \ell}u
+
\nabla_j\nabla_{\bar p}\nabla_{\bar k}u
+
X_{\bar p}\nabla_j\nabla_{\bar k}u\bigg)
\nabla_q u \, g^{j\bar p}g^{q\bar k}\,e^{\theta}\o^n
\nonumber\\
&=&
\|\bar\nabla\bar\nabla u\|_\theta^2
+
\int_M R_{\bar kj}\nabla^ju\nabla^{\bar k}u\,e^{\theta}\o^n
\nonumber\\
&&
\mbox{}  +
\int_M
X_j \nabla_{\bar p}\nabla_{\bar k}u
\,\nabla_qu \, g^{j\bar p}g^{q\bar k}
e^{\theta}\o^n
-
\int_M
X_{\bar p}\nabla_j\nabla_{\bar k}u
\,\nabla_qu \, g^{j\bar p}g^{q\bar k}
e^{\theta}\o^n. \ \ \ 
\eea
Now we rewrite the integrands of the last two terms in the last line
as follows:
\bea
X_j \nabla_{\bar p}\nabla_{\bar k}u
\,\nabla_qu g^{j\bar p}g^{q\bar k}
&=&
\nabla_{\bar k}(X^{\bar p}\nabla_{\bar p}u)
\nabla_q u g^{q\bar k}
-
(\nabla_{\bar k}X^{\bar p})\nabla_{\bar p}u\nabla_q u g^{q\bar k}
\nonumber\\
X_{\bar p}\nabla_j\nabla_{\bar k}u
\,\nabla_qu g^{j\bar p}g^{q\bar k}
&=&
\nabla_{\bar k}(X^j\nabla_ju)\nabla_q u g^{q\bar k},
\eea
where we have used the fact that $X^j$ is a holomorphic vector field.
But
\be
X^{\bar p}\nabla_{\bar p}u-X^j\nabla_ju
=
\bar X u-Xu=0,
\ee
since $u$ is invariant under ${\rm Im}\,X$ and thus 
we are left with the desired formula (\ref{bk}).
Substituting this formula and the relation (\ref{zero}) in the earlier
identity (\ref{identity}) gives the identity (\ref{analogue1}).

\medskip
Once the identity (\ref{analogue1}) is available,
the arguments of \cite{PSSW1} apply to give the proof of Theorem
\ref{differencedifferential}, with suitable modifications.
Write $\pi(\overline{\nabla} u)$ for the orthogonal projection with respect to the norm $ \| \cdot \|_{\theta}$  
of the $T^{1,0}$ vector field $\bar \nabla u$ onto the space of holomorphic vector fields.
Then
\bea
\label{old49}
\|\bar\nabla\bar\nabla u\|_\theta^2
&\geq&
\lambda_X(t)\|\bar\nabla u-\pi(\bar\nabla u)\|_\theta^2
\nonumber\\
&=&
\lambda_X(t)\,\bigg(\|\bar\nabla u\|_\theta^2  - \|\pi(\bar\nabla u)\|_\theta^2\bigg),
\eea
where $\lambda_X(t)$ is the eigenvalue introduced in (\ref{lambdaX}).
Since 
$\|\pi(\bar \nabla u)\|_{\theta}^2=\int_M \pi(\bar\nabla u)^j\pl_ju \, e^{\theta}
\omega^n= -F_X(\pi\bar\nabla u)$,
we obtain the inequality
\bea 
\label{YXidentity1}\nonumber
\dot Y_X(t)
&\leq&
-2 \lambda_X(t)\,Y_X(t)- 2\lambda_X(t)\,F_X(\pi\bar\nabla u)
+
\int_M|\na u|^2(Xu)e^\theta\o^n
\\
&&-
\int_M(R_{\bar kj}-g_{\bar kj}-\na_jX_{\bar k})\na^j u\overline{\na u^k}\,e^\theta\o^n
-
\int_M (R-n-\na_jX^j)|\na u|^2 e^\theta\o^n. \ \ \  \ \ 
\eea

We come now to the proof of the difference-differential inequality for $Y_X(t)$
in the statement of Theorem \ref{differencedifferential}.
First, observe that
\be
\|R_{\bar kj}-g_{\bar kj}-\na_jX_{\bar k}\|_{L^2}
=
\|R-n-\na_jX^j\|_{L^2}.
\ee
This is because the left hand side equals $\|\na\bar\na u\|_{L^2}$
and the right hand side equals $\|\Delta u\|_{L^2}$, which are readily
seen to be equal by an integration by parts. Next, we claim
that the last three terms on the right hand side of (\ref{YXidentity1})
can all be bounded by a constant multiple of
\be
\,\|\na u\|_{L^2}\,\|(u-b)(t-2)\|_{C^0}^2.
\ee
Indeed, since $\theta$ is bounded, we can write
\bea
\left| \int_M(R_{\bar kj}-g_{\bar kj}-\na_jX_{\bar k})\na^ju\overline{\na^k u}
e^\theta\o^n \right|
&\leq& C \|\na u\|_{C^0}\|\na u\|_{L^2}\|R_{\bar kj}-g_{\bar kj}-\na_jX_{\bar k}\|_{L^2}
\nonumber\\
&=& C \|\na u\|_{C^0}\|\na u\|_{L^2}\|R-n-\na_jX^j\|_{L^2}
\nonumber\\ \label{term1}
&\leq&
C\,\|\na u\|_{L^2}\,\|(u-b)(t-2)\|_{C^0}^2,
\eea
where we have applied Lemma \ref{smoothinglemma}  to obtain the last line.  Note that if $\| (u-b)(t-2) \|_{C^0}>\varepsilon$, for $\varepsilon$ as in Lemma \ref{smoothinglemma}, then we can still obtain the bound
\be
\|\na u\|_{C^0}\|R-n-\na_jX^j\|_{L^2} \le
C\, \|(u-b)(t-2)\|_{C^0}^2,
\ee
using the uniform estimates of $\| \nabla u\|_{C^0}$ and $\| \Delta u\|_{C^0}$.
Similarly,
\be \label{term2}
\left| \int_M(R-n-\na_jX^j)|\na u|^2e^\theta\o^n \right|
\leq
C\,\|\na u\|_{L^2}\,\|(u-b)(t-2)\|_{C^0}^2,
\ee
while the estimate for the remaining term,
\bea \label{term3}
\left| \int_M|\na u|^2(Xu)e^\theta\o^n \right|
\leq
\|\na u\|_{L^2}\,\|(u-b)(t-2)\|_{C^0}^2,
\eea
is even easier, since $|Xu|\leq |X|\cdot|\na u|
\leq C\,\|(u-b)(t-2)\|_{C^0}$.

\smallskip
Let $0<\rho := 1/( n+1)<1$. By the $L^2/C^0$ Poincar\'e inequality
and Lemma \ref{smoothinglemma}, we can write
\bea
\|(u-b)(t-2)\|_{C^0}^2
&\leq& C\,\|\na u(t-2)\|_{L^2}^{2\rho}
\|\na u(t-2)\|_{C^0}^{2(1-\rho)}
\nonumber\\ 
&\leq&
C \, Y_X(t-2)^{\rho}\|(u-b)(t-4)\|_{C^0}^{2(1-\rho)}.
\eea
We note that these inequalities are homogeneous, in the sense that the sum
of the exponents on either side always match. We can thus iterate, and
obtain
\bea \nonumber
\lefteqn{
\|(u-b)(t-2)\|_{C^0}^2
} \\
& \leq &  C\, Y_X (t-2)^\rho
\|(u-b)(t-4)\|^{2(1-\rho)}
\nonumber\\
&\leq&
C\,Y_X(t-2)^\rho Y_X(t-4)^{2(1-\rho)\rho}
\|(u-b)(t-6)\|_{C^0}^{2(1-\rho)^2}
\nonumber\\
&\leq& \cdots
\nonumber\\
&\leq&
C\, Y_X(t-2)^{\frac{\delta_1}{2}}Y_X(t-4)^{\frac{\delta_2}{2}}
\cdots Y_X(t-2N)^{\frac{\delta_N}{2}}
\|(u-b)(t-2(N+1))\|_{C^0}^{2(1-\rho)^N}, \ \ \ \ \label{iterate}
\eea
with $\sum_{j=1}^N\delta_j+2(1-\rho)^N=2$.
Fix $N$ with $2(1-\rho)^N<1$ and set $\delta_0=1$.  Since the quantity $\|(u-b)(t-2(N+1))\|_{C^0}$
is bounded by Lemma \ref{lemmaddtphi}, the statement of Theorem \ref{differencedifferential} follows from the inequalities (\ref{YXidentity1}), (\ref{term1}), (\ref{term2}), (\ref{term3}) and (\ref{iterate}).

\section{Proof of Theorem \ref{Rconvergence}} \label{PT2}
\setcounter{equation}{0}


Theorem \ref{Rconvergence} follows easily from what we have proved above.
Indeed, by (\ref{Kenergy1}) and   (\ref{phiu}),
the variation of the modified Mabuchi energy along the modified
K\"ahler-Ricci flow is given by
\be
\dot\mu_X=
-{1\over V}\int_M |\na u_{X,\o}|^2\,e^{\theta_{X,\o}}\o^n
=
-{1\over V}Y_X(t).
\ee
Integrating in $t$, we see that
 condition ($A_X$) implies:
\be \label{Y1}
\int_0^{\infty} Y_X (t) dt < \infty.
\ee
On the other hand, from (\ref{identity}) and
the uniform bounds of  $\theta$, $X u_{X,\o}$, $R$ and $\nabla_j X^j$ we obtain
\be \label{Y2}
\dot Y_X \le C Y_X.
\ee
The inequalities (\ref{Y1}) and (\ref{Y2}) imply (as in Section \S2 of \cite{PS1}) that $Y_X(t) \rightarrow 0$ as $t \rightarrow \infty$.

By the uniform bound of $\| \nabla u_{X,\o}\|_{C^0}$ and  Lemma \ref{l3} we have
\be
\| u_{X,\o} - b \|_{C^0} \rightarrow 0, \quad \textrm{as } t \rightarrow \infty.
\ee
Then from Lemma \ref{smoothinglemma} we see that
\be
\| \Delta u_{X,\o} \|_{C^0} \rightarrow 0, \quad \textrm{as } t \rightarrow \infty.
\ee
Since $\Delta u_{X,\o} = R - n - \nabla_j X^j$,
the first part of Theorem \ref{Rconvergence} is established.
The $L^p$ integrability of $\|R-n-\na_jX^j\|_{C^0}$ on $[0,\infty)$
is established in the same way as part (ii) of Theorem 1
in \cite{PSSW1}. The proof of Theorem \ref{Rconvergence}
is complete.

\section{Proof of Theorem \ref{conditions}} \label{PT1}
\setcounter{equation}{0}

It is convenient to introduce the following fifth condition:

\smallskip
(o) For each $k =0, 1,2, \ldots $, there exists a finite constant $A_k$ so that
\be
{\rm sup}_{t\geq 0}\,\|\varphi\|_{C^k}
\ \leq\ A_k.
\ee

\smallskip
\noindent
We shall prove the following implications
\bea
\label{implications}
&&
(o)\Leftrightarrow (iii)
\nonumber\\
&&
(o) \Rightarrow (iv) \Rightarrow (ii) \Rightarrow (iii)
\nonumber\\
&&
(iv)\Rightarrow (v)\Rightarrow (i)\Rightarrow (iii)
\eea
from which the equivalence of all five conditions
(i)-(v), and hence Theorem \ref{conditions}, all follow at once.

\subsection{(o) $\Leftrightarrow$ (iii)}

This is the extension to the case of the modified K\"ahler-Ricci flow
of the classical fact that a $C^0$ estimate for the complex Monge-Amp\`ere equation
implies $C^k$ estimates to all orders. We present it here, with emphasis
only on those points that require additional arguments.  We note that in \cite{TZ2}, a different method is used to obtain higher order estimates, involving a modification of the potential $\varphi$ along the flow.  We give here a direct proof of the higher order estimates for solutions of (\ref{KRFP2}).

\medskip

The first step is to show that $C^0$ estimates for $\varphi$ imply second order
estimates for $\varphi$.  In this section, for ease of notation, we will use $\hat{g}_{\bar k j}$
to denote the original metric $g^0_{\bar k j}$, and  $\hat{\Delta}$ for the Laplacian with
respect to this metric.
As in the original approach of Yau \cite{Y1} and
Aubin \cite{A}, we apply the maximum principle to the flow
of $\log (n+\hat\Delta\varphi)- A\varphi$, where $A$ is a large constant to
be chosen later. It is convenient to use the formulas obtained
in \cite{PSS} for general flows. For this, we introduce the endomorphism
\be
h^\al{}_\b=\hat g^{\al\bar\g}g_{\bar\g\b}.
\ee
Then $n+\hat\Delta\varphi={\rm Tr}\,h$, and we have
(see e.g. \cite{PSS}, eq. (2.27))
\bea
(\Delta-\pl_t)\log\,{\rm Tr}\,h
&=&
{1\over {\rm Tr}\,h}
\bigg(\hat\Delta(\log{\o^n\over\o_0^n}-\dot\varphi)-\hat R\bigg)
\nonumber\\
&&
\mbox{} +
{1\over {\rm Tr}\,h}
h^r{}_j(h^{-1})^p{}_s\hat g^{s\bar q}\hat R_{\bar qp}{}^j{}_r
\nonumber\\
&&
\mbox{} +
\bigg\{
{g^{j\bar k}{\rm Tr}(\na_j h\,h^{-1}\,\na_{\bar k}h)
\over {\rm Tr}\,h}
-
{g^{j\bar k}\na_j{\rm Tr}\,h
\na_{\bar k}{\rm Tr}\,h\over ({\rm Tr}\,h)^2}\bigg\},
\eea
where $\hat R$, $\hat R_{\bar kj}{}^\al{}_\b$ are the scalar and the Riemann
curvature tensor of the metric $\hat g_{\bar kj}$, $\hat g^{i\bar k}$ is its
inverse, and otherwise all indices are raised and lowered using the metric $g_{\bar kj}$.
The last line is non-negative (see \cite{Y1}), and the second
line is bounded below by $-C_1{\rm Tr}\,h^{-1}$. Since
$0<({\rm Tr}\,h)^{-1}<{\rm Tr}\,h^{-1}$, we obtain
\bea
(\Delta-\pl_t)\log\,{\rm Tr}\,h
&\geq&
{1\over {\rm Tr}\,h}
\hat\Delta(\log{\o^n\over\o_0^n}-\dot\varphi)
-C_1\,{\rm Tr}\,h^{-1}.
\eea
For the modified K\"ahler-Ricci flow, we have
\be
\hat\Delta(\log{\o^n\over\o_0^n}-\dot\varphi)
=
-\hat\Delta\varphi-\hat\Delta\theta - \hat{\Delta} f(\omega_0)
=
-{\rm Tr}\,h+n-\hat\Delta\theta - \hat{\Delta} f(\omega_0).
\ee
The new term compared to the K\"ahler-Ricci flow is $-\hat\Delta\theta$,
which is not yet known to be bounded. To eliminate it, we consider
instead the expression
$(\Delta+X-\pl_t)\log\,{\rm Tr}\,h$.
The main observation is that:
\be
\label{XTrh}
X\,{\rm Tr}\,h
= X (\hat{\Delta} \varphi) =
\hat\Delta \theta-\hat\Delta \hat{\theta}
+(\hat\na_jX^m) \,h^j{}_m-\hat\na_mX^m,
\ee
where $\hat{\theta}= \theta_{X, \omega_0}$.
To see this, observe that using the fact that $X$ is holomorphic,
\bea
\hat{\Delta} (X \varphi) & = & \hat{g}^{j \bar k} \hat{\nabla}_j \hat{\nabla}_{\ov{k}}
(X^m \hat{\nabla}_m \varphi)
 =  g^{j \bar k} \hat{\nabla}_j X^m (g_{\bar k m} - \hat{g}_{\bar k m}) + X( \hat{\Delta} \varphi).
\eea
Then (\ref{XTrh}) follows from the identity $X \varphi = \theta-\hat{\theta}$.

Hence for the modified K\"ahler-Ricci flow,
\bea
(\Delta+X-\pl_t)\log\,{\rm Tr}\,h
&\geq&
{1\over {\rm Tr}\,h}
(-{\rm Tr}\,h+n-\hat\Delta\hat{\theta}
+(\hat\na_jX^m) \,h^j{}_m-\hat\na_mX^m)
-C_1\,{\rm Tr}\,h^{-1}
\nonumber\\
&\geq&
-C_2-C_3\,{\rm Tr}\,h^{-1}.
\eea
From here, the proof can proceed as before. Set $A=C_3+1$. Since $\Delta\varphi=-{\rm Tr}\,h^{-1}+n$, and
$\dot\varphi$ and $X\varphi=\theta-\hat{\theta}$ are bounded by Lemmas \ref{lemmaddtphi}
and \ref{lemmathetaX}, we have
\bea \label{secondorder}
(\Delta+X-\pl_t)(\log\,{\rm Tr}\,h-A\varphi)
&\geq&
-C_4+\,{\rm Tr}\,h^{-1}.
\eea
 Fix $T>0$.  Then at a maximum point $(x_0, t_0)$ of the function
$\log\,{\rm Tr}\,h-A\varphi$ on $M \times (0,T]$, the quantity ${\rm Tr}\, {h^{-1}}$ is bounded from above.  But from the modified K\"ahler-Ricci
flow equation (\ref{KRFP2}) and the fact that $\dot\varphi$, $\theta$ and $\varphi$ are bounded, we see that the logarithm of the product of the eigenvalues of $h$ is bounded, giving an upper bound of ${\rm Tr} \, h(x_0, t_0)$.  Thus, using again the estimate of $\| \varphi\|_{C^0}$, we obtain a uniform upper bound of ${\rm Tr}\, h$ along the modified K\"ahler-Ricci flow.  The positivity of the metric  $g_{\bar k j} = \hat{g}_{\bar k j} + \partial_j \partial_{\bar k} \varphi$ ensures uniform bounds of $\partial_j \partial_{\bar k} \varphi$ along the flow.  In addition, from the lower bound of $\log ( \omega^n/ \hat{\o}^n)$, we have an estimate $g_{\bar k j} \ge C_5 \hat{g}_{\bar k j}$, for $C_5>0$, showing that  $g$ is uniformly equivalent to $\hat{g}$ along the flow.

\medskip

We now give the third order estimates.
As in \cite{Y1}, set
$\varphi_{j\bar k m}=\hna_m\pl_{\bar k}\pl_j\varphi$ and
$S\equiv
g^{j\bar r} g^{s\bar k}g^{m\bar t}  \varphi_{j\bar k m}\varphi_{\bar r s\bar t}$.
Again, it is convenient to compute instead in terms of the
connection $\na h h^{-1}$, in terms of which we have
\be
S=
g^{m\bar\g}g_{\bar\mu \b}g^{\ell\bar\al}
(\nabla_mh\,h^{-1})^\b{}_{\ell}\overline{(\nabla_\g h\,h^{-1})^\mu{}_\al}
=
|\nabla h\,h^{-1}|^2.
\ee
From \cite{PSS} eq. (2.48), under any
flow, we have the general formula
\bea
(\Delta-\pl_t) S&=& |\bar\nabla(\nabla h\,h^{-1})|^2+|\nabla(\nabla h\,h^{-1})|^2
\nonumber\\
&&
\mbox{} +
g^{m\bar\g}\<(\Delta-\pl_t)(\nabla_m h\,h^{-1}),\nabla_\g h\,h^{-1}\>
+
g^{m\bar\g}\<\nabla_m h\,h^{-1},(\Delta-\pl_t)(\nabla_\g h\,h^{-1})\>
\nonumber\\
&&
\mbox{} +
\bigg((h^{-1}\dot h+{\rm Ric})^{m\bar \g}g_{\bar\mu\beta}g^{\ell\bar\alpha}
-
g^{m\bar\g}(h^{-1}\dot h+{\rm Ric})_{\bar\mu \b}g^{\ell \bar\alpha}
+
g^{m\bar\g}g_{\bar\mu\b}(h^{-1}\dot h+{\rm Ric})^{\ell\bar\alpha}\bigg)
\nonumber\\
&&
\qquad\qquad\times
(\nabla h_m\,h^{-1})^\b{}_{ \ell} \overline{(\nabla_\g h\,h^{-1})^\mu{}_\alpha}.
\eea
Next, we specialize to the modified K\"ahler-Ricci flow.
We always have the following formula relating the curvatures
of two metrics $g_{\bar kj}$ and $\hat g_{\bar kj}$
\be
\Delta(\nabla_m h\,h^{-1})^\alpha{}_\beta
=
\nabla^{\bar q}\hat R_{\bar qm}{}^\alpha{}_\b-\nabla_m R^\alpha{}_\beta.
\ee
For the modified K\"ahler-Ricci flow, we have
\be
\pl_t(\nabla_m h\,h^{-1})^\alpha{}_\beta
=
\nabla_m(h^{-1}\dot h)^\alpha{}_\beta
=
-
\nabla_m R^\alpha{}_\beta
+
\nabla_m\nabla_\b X^\alpha.
\ee
Thus we obtain
\be
\label{modifiedcalabi}
(\Delta -\pl_t)(\nabla_m h \,h^{-1})^\alpha{}_\b
=
\nabla^{\bar q}\hat R_{\bar qm}{}^\alpha{}_\b
-
\nabla_m\nabla_\b X^\alpha.
\ee
Similarly
\be
(h^{-1}\dot h)_{\bar kj}+R_{\bar kj}=g_{\bar kj}+\na_j X_{\bar k}.
\ee
Clearly, the terms $\nabla_m\nabla_\b X^\al$ and
$\na_jX_{\bar k}$ on the right hand side of the previous two equations
are the only changes due to the modified K\"ahler-Ricci flow.
Using then the equation (2.51) for the K\"ahler-Ricci flow in \cite{PSS},
we obtain immediately the following formula
\bea
(\Delta-\pl_t)S
&=&
|\bar\nabla(\nabla h\,h^{-1})|^2+|\nabla(\nabla h\,h^{-1})|^2
+|\nabla h\,h^{-1}|^2
\nonumber\\
&&
\mbox{} +
g^{m\bar\g}\nabla^{\bar q}\hat R_{\bar q m}{}^\b{}_{\ell}
\overline{(\nabla_\g h\,h^{-1})_{\bar\b}{}^{\bar \ell}}
+
g^{m\bar \g}(\nabla_m h\,h^{-1})_{\bar\mu}{}^{\bar\al}
\overline{\nabla^{\bar q}\hat R_{\bar q\g}{}^\mu{}_\al}
\nonumber\\
&&
\mbox{} +
(I)+(II)+(III)+(IV)+(V),
\eea
where the terms (I)-(V) are due to the modifications arising from the holomorphic vector field $X$,
and given explicitly by
\bea
(I)& = &
\na^{\bar\g}X^m
g_{\bar\mu \b}g^{\ell\bar\al}
(\na_m h\,h^{-1})^\b{}_\ell
\overline{(\na_\g h\,h^{-1})^\mu{}_\al}
\nonumber\\
(II) &= &
-g^{m\bar\g}g_{\bar\mu\b}g^{\ell\bar\al}\na_m\na_\ell X^\b
\overline{(\na_\g h\,h^{-1})^\mu{}_\al}
\nonumber\\
(III) & = &
g^{m\bar \g}g_{\bar\mu\b}
\na^{\bar\al}X^\ell
(\na_m h\,h^{-1})^\b{}_\ell
\overline{(\na_\g h\,h^{-1})^\mu{}_\al}
\nonumber\\
(IV) & =&
-g^{m\bar\g}g_{\bar\mu\beta}g^{\ell\bar\al}
(\na_m h\,h^{-1})^\b{}_\ell
\overline{\na_\g\na_\al X^\mu}
\nonumber\\
(V) & = &
-g^{m\bar\g}\na_\b X_{\bar\mu}
g^{\ell\bar\al}
(\na_m h\,h^{-1})^\b{}_\ell
\overline{(\na_\g h\,h^{-1})^\mu{}_\al}.
\eea
Because of the presence of the connection in e.g. $\na^{\bar\g}X^m
=g^{\al\bar \g}\na_\al X^m$, the first covariant derivatives of
$X^m$ are of order $O(S^{1\over 2})$, and hence
\be
|(I)|+|(III)|+|(V)|
\leq C_{6}\,S\,|\na X|.
\ee
The second covariant derivatives
of $X^m$ can be expressed as follows
\bea
\na_\g\na_\al X^\mu
&=&
\na_\g(\hna_\al X^m+(\na_\al h\,h^{-1})^\mu{}_\nu X^\nu)
\nonumber\\
&=&
\na_\g\hna_\al X^m
+
\na_\g(\na_\al h\,h^{-1})^\mu{}_\nu X^\nu+
(\na_\al h\,h^{-1})^\mu{}_\nu\na_\g X^\nu.
\eea
The first term on the right hand side is again of order $O(S^{1\over 2})$.
The second term can be bounded by $|\na(\na h\,h^{-1})|$ since $|X|$
is bounded. Thus we can write
\footnote{There are actually some partial cancellations between the
terms I-V. We shall not need this fact, and won't exhibit it more
explicitly.}
\bea
|(II)|+|(IV)|
&\leq&
C_{7}\, (S+1) +|X|\cdot|\na(\na h\,h^{-1})|\cdot|\na h\,h^{-1}|
+
|\na h\,h^{-1}|^2\cdot|\na X|
\nonumber\\
&\leq&
C_{8}\,S+{1\over 2}|\na(\na h\,h^{-1})|^2+S\,|\na X|+ C_{9}.
\eea
Putting this all together, we obtain the following estimate for
the flow of $S$ in the modified K\"ahler-Ricci flow,
\bea \label{evolveS}
(\Delta-\pl_t)S
\geq
{1\over 2}|\na(\na h\,h^{-1})|^2
+
|\bar\na(\na h\,h^{-1})|^2
-
C_{9}\,S\,|\na X|
-
C_{10}(1+S).
\eea

By the method of \cite{Y1}, we can control terms of order $O(S)$ using the evolution equation for
$\textrm{Tr}\, h$.  However, we will need an additional argument to deal with the quantity
 $S | \nabla X|$ which is of the order $O(S^{3/2})$.   Since $|X|$ is uniformly bounded along the flow,
we have
 \begin{eqnarray} \label{evolveX}
(\Delta - \partial_t )   |X|^2 & = &  | \nabla X|^2 - |X|^2 - \partial_i \partial_{\ov{j}}
\theta X^i \ov{X^j}
 \ge    \frac{1}{2} | \nabla X|^2 - C_{11}.
\end{eqnarray}

We define a constant $K= 65 \sup_{M \times [0, \infty)} ( |X|^2 +1)$ and compute the
evolution of the quantity $S/(K-|X|^2)$.  Combining (\ref{evolveS}) and (\ref{evolveX})
we have
\begin{eqnarray} \nonumber
(\Delta - \partial_t)   \left( \frac{S}{K- |X|^2} \right) & \ge &
\frac{ \left( | \nabla (\nabla h \cdot h^{-1}) |^2 - | \overline{\nabla} (\nabla h \cdot h^{-1}) |^2
\right)}{2(K-|X|^2)} + \frac{S | \nabla X|^2}{2 (K - |X|^2)^2} \\ \nonumber
&& \mbox{}  + \frac{2 \Re (g^{i \ov{j}} \partial_i   S \, \partial_{\overline{j}} |X|^2)}
{(K-|X|^2)^2} + \frac{2S | \nabla |X|^2|^2}{(K-|X|^2)^3} \\
&& \mbox{} - \frac{C_{9} S | \nabla X|}{K- |X|^2} - C_{12}(1+S).
\end{eqnarray}
We will use the good first and second terms on the right hand side of this inequality
to deal with the bad third and fifth terms.  We estimate the third term as follows:
\begin{eqnarray} \nonumber
\frac{|2  g^{i \ov{j}} \partial_i   S \, \partial_{\overline{j}} |X|^2 |}{(K-|X|^2)^2} & \le &
\frac{ S | \nabla |X| \, |^2}{4(K-|X|^2)^2} + \frac{32 |X|^2 \left( | \nabla (\nabla h \cdot h^{-1}) |^2
+ | \overline{\nabla} (\nabla h \cdot h^{-1}) |^2  \right)}{(K- | X|^2)^2} \\
& \le & \frac{ S | \nabla X |^2}{4(K-|X|^2)^2} + \frac{ \left( | \nabla (\nabla h \cdot h^{-1}) |^2 +
| \overline{\nabla} (\nabla h \cdot h^{-1}) |^2  \right)}{2(K- | X|^2)}.
\end{eqnarray}
For the fifth term, observe that:
\begin{equation}
\frac{C_{9} S | \nabla X|}{K- |X|^2} \le \frac{S | \nabla X|^2}{4(K-|X|^2)^2} + C_{9}^2 S.
\end{equation}
Combining all of the above, we obtain
\begin{equation}
(\Delta -  \partial_t)  \left( \frac{S}{K- |X|^2} \right) \ge - C_{13}(1+S).
\end{equation}
But from the computation for the second order estimate, we have
\be
(\Delta - \partial_t) \textrm{Tr} \, h \ge  \frac{1}{2} S - C_{14},
\ee
and so applying the maximum principle to the quantity $\displaystyle{\left( S/(K-|X|^2) + 3C_{13} \,
\textrm{Tr} \, h \right)}$ it follows that $S/(K-|X|^2)$ and hence $S$ is bounded.

\bigskip
\noindent
{\bf Remark} \ Instead of computing the evolution of $S/(K- |X|^2)$, an alternative is to compute the evolution of the tensor $T_j^k = (\nabla_j h \, h^{-1})^k_l X^l$.  Indeed one can prove that:
\be
(\Delta - \partial_t) | T|^2 \ge - B_1 S - B_2,
\ee
for uniform constants $B_1$ and $B_2$.  Combining this with the evolution of $\textrm{Tr} \, h$ gives an upper bound of $|T|$ and hence $|\nabla X|$ along the flow.  This implies that the term $S |\nabla X|$ is of order $O(S)$ and one can proceed in the usual way to bound $S$.

\bigskip

In order to apply the standard parabolic estimates to obtain the higher order estimates, we require a derivative bound of $g_{\bar k j}$ in the $t$-direction (cf. \cite{Ch}, for example).  Given the estimates proved so far, it is sufficient to bound $| \textrm{Ric}(g) |$.  The evolution of the Ricci curvature along the modified K\"ahler-Ricci flow is given by
\be
( g^{p \ov{q}} \nabla_p \nabla_{\ov{q}} - \partial_t) R_{\ov{k} j} = R_{\ov{k}}^{ \ \ \ov{\ell}} R_{\ov{\ell} j} - R_{\ov{k} j}^{\ \ \ p \ov{q}} R_{\ov{q} p} - X^{\ell} \nabla_{\ell} R_{\ov{k} j} + R_{\ov{k} \ell } \nabla_j X^{\ell}.
\ee
Then given the estimates on $\textrm{Tr} \, h$ and $S$, we have
\bea \nonumber
( \Delta -\partial_t + X) | \textrm{Ric}(g) | & =&  \frac{1}{| \textrm{Ric} (g) |} \left\{ | \nabla \textrm{Ric}(g) |^2 - | \nabla | \textrm{Ric} | \, |^2 + | \textrm{Ric}(g) |^2  \right. \\ \nonumber
&& \left. \mbox{} - R_{\ov{k} j}^{\ \ \ r \ov{s}} R_{\ov{s} r} R^{j \ov{k}} + \nabla^{\ov{k}} X^p g^{j \ov{q}} R_{\ov{k} j} R_{\ov{q} p} - R_{\ov{k} \ell} \nabla_j X^{\ell} R^{j \ov{k}} \right\} \\
& \ge & - C_{15} ( | \textrm{Rm} |^2 +1).
\eea
But from the computations above for the evolution of $S$, there exist uniform constants $C_{16}$, $C_{17}$ with $C_{16}>0$  
such that
\be
( \Delta - \partial_t + X) S \ge C_{16} | \textrm{Rm} |^2 -C_{17}.
\ee
Then by applying the maximum principle to the quantity $| \textrm{Ric} (g) | + \frac{1}{C_{16}} (C_{15}+1) S$ we obtain the desired upper bound on $| \textrm{Ric}(g) |$.

\medskip

We have now established uniform parabolic $C^1$ estimates for the metric $g_{\ov{k} j}$ along the flow.
One can now obtain the higher order estimates in the usual way.   We differentiate the equation (\ref{KRFP2}) in space, making use of Lemma \ref{lemmaTZ}, and then apply the standard parabolic estimates (se \cite{L}, for example) together with a bootstrapping argument.

\subsection{(o) $\Rightarrow$ (iv)}

The remaining implications are all straightforward adaptations of arguments
in \cite{PSSW1}, so we shall be brief.
It is convenient to formulate the following lemmas:

\begin{lemma}
\label{decay}
Let $W(t)$ be a non-negative $C^\infty$ function of $t\in[0,\infty)$ with $W(t) \le K_0$ 
satisfying the difference-differential inequality
\be
\dot W(t)
\leq -2 \lambda\,W(t)+\lambda \prod_{j=0}^NW(t-2j)^{\nu_j\over 2}, \ \ \textrm{for } t \ge K_1 \ge 2N,
\ee
where $\lambda$ is a strictly positive constant,
and $\nu_j\geq 0$ satisfy ${1\over 2}\sum_{j=0}^N\nu_j=1$. Then there
exist constants $C,\kappa$ with $\kappa>0$ depending only on $K_0, K_1, \lambda, N, \nu_j$ so that
\be
W(t)\leq C\,e^{-\kappa t}.
\ee
\end{lemma}

\noindent
{\it Proof of Lemma \ref{decay}}: See Section \S5 of \cite{PSSW1}. Q.E.D.

\begin{lemma} \label{equivalencelemma} 
There exist constants $c_1, c_2>0$ depending only on the complex manifold $M$ and the holomorphic vector field $X$ such that for all $\omega \in {\cal K}_X$,
\be
c_1 \lambda(\omega) \le \lambda_X(\omega) \le c_2 \lambda(\omega).
\ee
\end{lemma}

\noindent
{\it Proof of Lemma \ref{equivalencelemma}}: This follows immediately from the fact that $\theta_{X, \o}$ is uniformly bounded (by a constant independent of choice of metric in ${\cal K}_X$) and the  argument of Lemma 1 in \cite{PSSW2}. Q.E.D.

\begin{lemma}
\label{formulationlemma}
Let $Y_X(t)$ be given as in (\ref{YX}) for the modified K\"ahler-Ricci
flow. Assume that the following three conditions
are satisfied

{\rm (a)} $F_X\equiv 0$

{\rm (b)} $Y_X(t)\to 0$ as $t\to\infty$

{\rm (c)} ${\rm inf}_{t\geq 0}\lambda(t) >0$.

\noindent
Then there exists constants $C,\kappa$ with $\kappa>0$ so that $Y_X(t)\leq C e^{-\kappa t}$.
\end{lemma}

\noindent
{\it Proof of Lemma \ref{formulationlemma}:} 
We apply
Theorem \ref{differencedifferential}. Under conditions (a)-(c),
the difference-differential inequality given there, together with Lemma \ref{equivalencelemma},  implies
that $Y_X(t)$ satisfies a difference-differential inequality
exactly of the type formulated in Lemma \ref{decay}.
The desired inequality follows then from this lemma. Q.E.D.

\medskip

Returning to the proof of (o) $\Rightarrow$ (iv),
assume that (o) holds, that is,
all norms $\|\varphi(t)\|_{C^k}$ are uniformly bounded in time for each $k$.
Then there exists a sequence of times $t_m \rightarrow \infty$ such that $\varphi(t_m)$ converges in $C^{\infty}$ to an element $\varphi(\infty)$ of ${\mathcal P}_X(M, \omega_0)$.  Since the modified Mabuchi K-energy $\mu_X$ is decreasing along the modified flow, it follows that for any time $t$,
\be
\mu_X (\varphi(t)) \ge \mu_X(\varphi(\infty)),
\ee
and hence $\mu_X$ is bounded below along the modified flow.  Using this, one can see that the limit metric $g_{\bar k j} (\infty)$ must be a K\"ahler-Ricci soliton with respect to $X$ (c.f. the proof of Theorem \ref{Rconvergence}).  Theorem \ref{theoremlowerbound} then establishes condition $(A_X)$. Next,
we claim that Condition (S)
is also satisfied, that is, the eigenvalues $\lambda(t)$ are bounded below by a strictly
positive constant $\lambda$. Otherwise, let $\varphi(t_m)$ be a subsequence with $\lambda(t_m)\to 0$.
It contains a subsequence $\varphi(t_{\ell})$ such that the corresponding metrics $g_{\bar kj}(t_{\ell})$
converge in $C^\infty$ to a K\"ahler-Ricci soliton
$g_{\bar kj}(\infty)$ with respect to $X$. In \cite{PS1}, it was shown that $\lambda(t_{\ell})\to\lambda(\infty)$
if $g_{\bar kj}(t_{\ell})\to g_{\bar kj}(t_\infty)$ and the dimensions of the holomorphic
vector fields of the complex structures for $g_{\bar kj}(t_{\ell})$ and $g_{\bar kj}(\infty)$
are the same. In the present case, the complex structures of $g_{\bar kj}(t_{\ell})$ and
$g_{\bar kj}(\infty)$ are the same, so we do have $\lambda(t_{\ell})\to\lambda(\infty)$.
Since $\lambda(\infty)>0$ by definition, we obtain a contradiction.
Condition $(S)$ is established.

\smallskip

The existence of a K\"ahler-Ricci soliton with respect to $X$ implies that (a) in
  Lemma \ref{formulationlemma} holds and condition $(A_X)$ gives (b) by Theorem \ref{Rconvergence}.
Since (c) in this Lemma is the same as $(S)$,
Lemma \ref{formulationlemma} applies,
and (iv) is established.

\subsection{(iv) $\Rightarrow$ (ii)}

Assume that (iv) is satisfied, and thus $Y_X(t)$ is rapidly decreasing.
Then Lemma \ref{l3} implies that
\be
\| u_{X,\o}-b \|_{C^0}\leq C\,e^{-{1\over 2(n+1)}\kappa t}
\ee
But Lemma \ref{smoothinglemma} implies then
that
\be
\|R-n-\na_jX^j\|_{C^0}\leq C'\, e^{-{1\over 2(n+1)}\kappa t}
\ee
which gives (ii).

\subsection{(ii) $\Rightarrow$ (iii)}

Assume (ii). Since $\pl_t\log (\o^n/\o_0^n)=g^{j\bar k}\dot g_{\bar kj}
=-(R-n-\na_jX^j)$,
we obtain immediately, uniformly in $t\in [0,\infty)$
\be
|\log ({\o^n\over\o_0^n})|\leq \int_0^\infty\|R-n-\na_jX^j\|_{C^0}\,dt
<C.
\ee
Next, from the modified K\"ahler-Ricci flow and the uniform bound for
$\|\dot\varphi\|_{C^0}$ (Lemma \ref{lemmaddtphi}), it follows that
\be
\|\varphi\|_{C^0}
=
\|\dot\varphi-\log({\o^n\over\o_0^n})-\theta+f(\omega_0)\|_{C^0}
\leq C.
\ee

\subsection{(iv) $\Rightarrow$ (v)}

Assume (iv). We have already seen that (iv) implies (ii), which implies
in turn (iii), which is equivalent to (o). Thus all metrics $g_{\bar kj}(t)$
are equivalent with uniform bounds in $t$. The same arguments as in \cite{PS1},
applied with straightforward modifications to $u_{X, \omega}=f+\theta_{X, \omega}$ instead of the Ricci potential
 there, show that $\|u_{X, \omega} \|_{(s)}$ converges exponentially fast
to $0$ with respect to any Sobolev norm $s$. It follows easily from there
that $g_{\bar kj}(t)$ converges exponentially fast to a K\"ahler-Ricci soliton
$g_{\bar kj}(\infty)$.

\bigskip
Since all the remaining implications in (\ref{implications}) are trivial,
the proof of Theorem \ref{conditions} is complete.

\medskip

\section{Proof of Theorem \ref{Fullconvergence}} \label{PT3}
\setcounter{equation}{0}

It is now easy to prove Theorem \ref{Fullconvergence}. If the modified K\"ahler-Ricci flow
converges to a K\"ahler-Ricci soliton with respect to $X$ then, by Theorem \ref{theoremlowerbound}, Condition
$(A_X)$ is satisfied. Furthermore, as part of the proof of the step
(o) $\Rightarrow$ (iv), we have seen that the uniform boundedness of
$\|\varphi(t)\|_{C^k}$ for each $k$ implies that Condition $(S)$
is satisfied. Thus it remains only to establish the sufficiency
of $(A_X)$ and $(S)$ for the exponential convergence of the K\"ahler-Ricci
flow.

\medskip
The arguments are now practically the same as in the proof of
the implication (o) $\Rightarrow$ (iv) earlier.  By Proposition \ref{Proposition} and Theorem \ref{Rconvergence}, $(A_X)$ implies (a) and (b) of  Lemma \ref{formulationlemma}.  In addition, $(S)$ gives
condition (c). Thus we obtain
the exponential decay of $Y_X(t)$, that is, Condition (iv)
of Theorem \ref{conditions} is satisfied. But Theorem \ref{conditions}
implies then the exponential convergence of the modified K\"ahler-Ricci
flow to a K\"ahler-Ricci soliton. Q.E.D.

\section{Further remarks and questions} \label{remarks}
\setcounter{equation}{0}

\noindent
\ (1) \ As mentioned in the introduction, our results imply that if a metric $\omega_0$ is invariant under the $S^1$ action induced by the imaginary part of a holomorphic vector field $X$ then
along the K\"ahler-Ricci flow starting at $\omega_0$, the quantity  $|X|^2$ is uniformly bounded.  We remark that an immediate consequence of this is that if $M$ and the initial metric are toric then $g_{\ov{k} j}(t)$ is bounded from above on compact subsets of the interior of the polytope $\Delta \subset {\mathbf R}^n$.  Of course, the case of the K\"ahler-Ricci flow on toric varieties is by now already well-understood by the work of Zhu \cite{Zhu2} (see also \cite{WZ}, \cite{TZ2}).

\bigskip
\noindent
\ (2) \   An obvious question is: what other notions of stability correspond to the existence of a K\"ahler-Ricci soliton on $M$?  By Theorem \ref{Rconvergence} we see that, assuming condition $(A_X)$, the $L^2$ norm of the quantity
$(R - n - \nabla_j X^j)$ tends to zero along the K\"ahler-Ricci flow.  Is there a natural analogue of the Calabi energy and a version of Donaldson's result \cite{D2} giving `semi-stability' in this case?

\bigskip
\noindent
\ (3) \ If one assumes uniform curvature bounds along the modified K\"ahler-Ricci flow then, in a similar vein to  \cite{PSSW2} and \cite{S2}, one would expect convergence of the flow to a soliton if the modified Futaki invariant $F_X$ vanishes and an addition stability assumption holds, such as: condition (S), condition (B) of \cite{PS1}, or a modified version of K-stability.

\bigskip
\noindent
\ (4) \ Does the stronger conclusion
\be 
\| R_{\bar k j} - g_{\bar k j} - \nabla_j X_{\bar k} \|_{C^0} \rightarrow 0, \ \ \textrm{as} \ t \rightarrow \infty
\ee
follow from condition $(A_X)$?  Observe that  $L^2$ convergence to zero does hold assuming $(A_X)$, since by Theorem \ref{Rconvergence},
\be
\int_M | \nabla \ov{\nabla} u_{X, \o} |^2 \omega^n = \int_M | \Delta u_{X, \o} |^2 \omega^n \rightarrow 0.
\ee

\bigskip
{}
\bigskip

\noindent
{\bf Acknowledgements}: The first author would like to thank the
Centre de Recerca Matematica in Barcelona and the Centre de Rencontres Mathematiques
in Luminy for their hospitality during the time when part of this research
was done.  In addition, the authors thank Gabor Szekelyhidi and Valentino Tosatti for some helpful discussions.

\bigskip

\pagebreak[3]

\v\v
$^*$ Department of Mathematics\\
Columbia University, New York, NY 10027\\

$^{**}$ Department of Mathematics \\
Rutgers University, Piscataway, NJ 08854\\

$^{\dagger}$ Department of Mathematics \\
Rutgers University, Newark, NJ 07102\\

$^{\ddagger}$ Department of Mathematics \\
University of California San Diego, La Jolla, CA 92093

\end{document}